\documentclass[final]{siamltex}
\usepackage{amsmath}
\usepackage{afterpage}
\usepackage[dvips]{graphicx}
\usepackage[normal]{subfigure}
\usepackage{amssymb}
\usepackage{subfig,float}
\usepackage{color}
\def\ni{\noindent}
\def\proof{{\ni \bf \underline{Proof:} }}
\def\endproof{\hfill$\blacksquare$\vspace{6pt}}

\newcommand{\bsub}{\begin{subequations}}
\newcommand{\esub}{\end{subequations}$\!$}

\newcommand{\ds}[0]{\displaystyle}

\newcommand{\R}{{\mathbb{R}}}
\newcommand{\eps}{{\displaystyle \varepsilon}}

\renewcommand{\theequation}{\arabic{section}.\arabic{equation}}

\newtheorem{example}{Example}[section]

\newcommand{\bex}{\begin{example}\rm}
\newcommand{\eex}{\end{example}}

\newcommand{\jl}[1]{\textcolor{black}{#1}}

\renewcommand{\theequation}{\arabic{section}.\arabic{equation}}
\renewcommand{\thesection}{\arabic{section}}

\title{Multiple Quenching Solutions of a Fourth Order Parabolic PDE with a singular nonlinearity modelling a MEMS Capac\jl{ito}r.}

\author{A. E. Lindsay \thanks{Department of Mathematics, University of Arizona, Tucson, Arizona, 85721, USA. ({\tt alindsay@math.arizona.edu})} \and J. Lega\thanks{Department of Mathematics, University of Arizona, Tucson, Arizona, 85721, USA. ({\tt lega@math.arizona.edu})} }

\begin{document}

\label{firstpage}
\maketitle

\begin{abstract}
Finite time singularity formation in a fourth order nonlinear parabolic \jl{partial differential equation (PDE)} is analyzed. The PDE is a variant of a ubiquitous model found in the field of Micro-Electro Mechanical System\jl{s} (MEMS) and is studied on a \jl{one-dimensional (1D)} strip and the unit disc. The solution itself remains continuous at the point of singularity while its higher derivatives diverge, a phenomenon known as quenching. For certain parameter regimes it is \jl{shown numerically} that the singularity will form at multiple isolated points in the 1D strip case and along a ring of points in the \jl{radially symmetric} 2D case. The location of these touchdown points \jl{is} accurately predicted by means of asymptotic expansions. The solution itself is shown to converge to a stable self-similar profile at the singularity point. Analytical calculations are verified by use of adaptive numerical methods which take advantage of symmetries exhibited by the underlying PDE to accurately resolve solutions very close to the singularity.
\end{abstract}



\pagestyle{myheadings}
\thispagestyle{plain}
\markboth{A. E. Lindsay \& J. Lega}{Multiple Touchdown in a MEMS Capacitor.}

\section{Introduction}

\begin{figure}[htbp]
\centering
\includegraphics[height=4cm,clip]{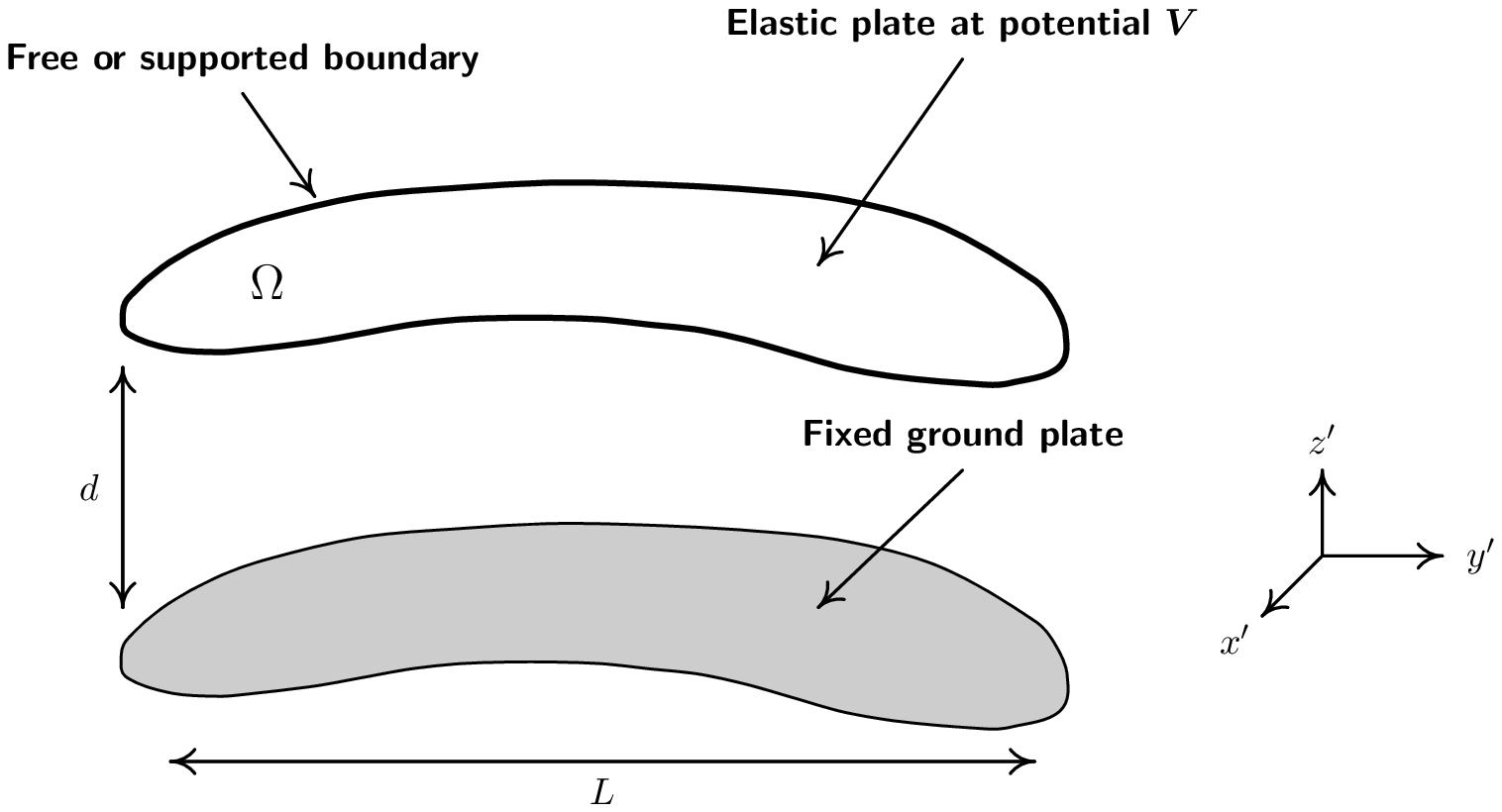} 
\parbox{4.5in}{\caption{Schematic plot of the MEMS capacitor (reproduced from \cite{LW}) with a deformable elastic
upper surface that deflects towards the fixed lower surface under an applied
voltage.\label{fig:schem}}}
\end{figure}

Micro-Electromechanical Systems (MEMS) combine electronics with
micro-size mechanical devices to design various types of microscopic
machinery (cf.~\cite{PB}). A key component of many MEMS is the
simple capacitor shown in Fig.~\ref{fig:schem}. The upper part of
this device consists of a thin deformable elastic plate that is held
clamped along its boundary, and which lies above a fixed ground plate.
When a voltage $V$ is applied between the plates, the upper surface can
exhibit a significant deflection towards the lower ground plate. When the applied voltage $V$ exceeds a critical value $V^{*}$, known as the \emph{pull-in} voltage, the deflecting surface can make contact with the ground plate. This phenomenon, known as \emph{touchdown}, will compromise the usefulness of some devices but is essential for the operation of others (\emph{e.g.} switches and valves). Capturing and quantifying this phenomenon is a topic of some mathematical interest and is the subject of this paper.

A canonical model, originally proposed in \cite{PB},
suggests the dimensionless deflection $u(x,t)$ of a device occupying a bounded region $\Omega\subset\mathbb{R}^2$ satisfies the fourth-order problem
\begin{equation}\label{bifull}
u_t = - \Delta^2 u + \delta \Delta u -
\frac{\lambda f(x)}{(1+u)^2} \,, \qquad x \in \Omega \,; \qquad \begin{array}{c} u = 0, \quad
\partial_n u = 0 \, \quad x \in \partial\Omega \,. \\[5pt] u = 0, \quad t=0 \quad x\in\Omega \end{array}
\end{equation}
Here, the positive constant $\delta$ represents the relative effects
of tension and rigidity on the deflecting plate, and $\lambda \geq 0$
represents the ratio of electric forces to elastic forces in the
system, and is directly proportional to the square of the voltage $V$
applied to the upper plate. The function $f(x)\in C^{\alpha}(\Omega)$ for $\alpha\in(0,1)$, represents possible heterogeneities in the deflecting surface's dielectric profile while the boundary conditions in (\ref{bifull})
assume that the upper plate is in a clamped state along its rim. The model (\ref{bifull}) was derived in \cite{PB} from a
narrow-gap asymptotic analysis.

The second order equivalent of \eqref{bifull} 
\begin{equation}\label{membrane_full}
u_t =  \Delta u -
\frac{\lambda f(x)}{(1+u)^2} \,, \qquad x \in \Omega \,; \qquad \begin{array}{c} u = 0, \quad x \in \partial\Omega \,. \\[5pt] u = 0, \quad t=0 \quad x\in\Omega \end{array}
\end{equation}
 has been the subject of extensive study recently and there are now many established results regarding the behaviour of solutions, both dynamic and steady (c.f. \cite{EGG2} and the references therein for a thorough account). In particular it is known that there exists a $\lambda^{\ast}>0$ such that whenever $\lambda>\lambda^{\ast}$ and $\inf_{\Omega}f>0$, the device touches down in finite time, \emph{i.e.} $||1+u(\cdot,t)||_{\inf} \to 0^{+}$ as $t\to t_c^{-}$. Lower and upper bounds have been established on the touchdown time $t_c$ of \eqref{membrane_full} and it is known that if touchdown occurs at an isolated $x_c\in\Omega$, then $f(x_c)\neq0$. Additionally, a refined asymptotic study of the touchdown profile was performed in \cite{GUO2} where it was shown that the quenching solution is not exactly self-similar and has asymptotic form
 \begin{equation}\label{intro_membrane_TD}
 u \to -1 + [3f(x_c)\lambda(t_c-t)]^{1/3}\left( 1 -\frac{1}{2|\log(t_c-t)|} + \frac{(x-x_c)^2}{4(t_c-t)|\log(t_c-t)|} +\cdots\right)
 \end{equation}
 where $x_c\in\Omega$ and $t_c>0$ are the touchdown location and time respectively. In addition, when $f(x)$ is a constant and $\Omega = [-1,1]$, the unique touchdown point is $x_c =0$.
 
In contrast to the second order problem \eqref{membrane_full}, very much less is know\jl{n} about the fourth order problem \eqref{bifull}, partly due to the lack of a maximum principle. In the absence of the tension term ($\delta = 0$) and with $f(x)=1$, equilibrium solutions of \eqref{bifull} were studied in \cite{CEGM} and the existence of a pull in voltage $\lambda^{\ast}$ was demonstrated for $\Omega$ a radially symmetric ball. The maximal branch of equilibrium solutions to \eqref{bifull}, \emph{i.e.} those solutions with largest $L_2$ norm for any sufficiently small $\lambda$, were constructed in the limit as $u\to-1^{+}$ in one and two dimensions in \cite{KLW,LW1}. Under the relaxed Navier boundary conditions $u=\Delta u = 0$ on $\partial\Omega$, a maximum principle is available and theoretical results regarding the existence and uniqueness of solutions are more tractable \cite{CEG,GW5}. 

Literature on the dynamics of fourth order MEMS equations is particularly sparse. The work of \cite{GUO} concerning the wave equation
\begin{equation}\label{GUOMAIN}
\begin{array}{cc}
 \mu w_{tt} + w_t - \Delta w + B \Delta^2 w = \frac{\lambda}{(1-w)^2}   & \mbox{in} \quad \Omega\times (0,T]\\[5pt]
w = \Delta w = 0 & \mbox{on} \quad \partial\Omega\times (0,T)\\[5pt]
w(x,0) = w_0(x), \quad w_t(x,0) = w_1(x) & \mbox{in} \quad \Omega
\end{array}
\end{equation}
appears to be the first contribution to the topic in which it is shown that \eqref{GUOMAIN} touches down in finite time for $\lambda>\lambda^{\ast}$.

In \jl{the present} work radially symmetric dynamical solutions of the fourth order MEMS problem
\bsub\label{mems_main}
\begin{equation}\label{mems_main_a}
u_t = -\eps^2 \Delta^2 u - \frac{1}{(1+u)^2}, \quad x\in\Omega, \qquad u(x,0) = 0, \quad x\in\Omega
\end{equation}
are considered for domains
\begin{equation}\label{mems_main_b}
\mbox{(Strip)}: \quad \Omega = [-1,1]; \qquad \mbox{(Unit Disc)}: \quad \Omega = \{ x^2 + y^2 \leq 1 \}
\end{equation}
and boundary conditions
\begin{equation}\label{mems_main_c}
\mbox{(Clamped)}: \quad u = 0, \quad \partial_n u = 0 \quad x\in\partial\Omega; \qquad \mbox{(Navier)}: \quad u = 0, \quad \Delta u = 0 \quad x\in \partial\Omega.
\end{equation}
\esub
The particular form of this equation is obtained from \eqref{bifull} by \jl{setting $f(x)=1$,} neglecting the tension term $\Delta u$ ($\delta = 0$), taking $\lambda t$ as a new time variable, and defining $\lambda=\eps^{-2}$. The consideration of radially symmetric solutions of \eqref{mems_main_a} on the strip and unit disc geometries effectively focusses attention on the PDE
\begin{equation}\label{mems_main_radial}
\jl{u_t = -\eps^2 \left[u'''' + \frac{2(N-1)}{r}u''' -\frac{N-1}{r^2} u'' + \frac{N-1}{r^3}u' \right] - \frac{1}{(1+u)^2}}
\end{equation}
for $N=1$ ( Strip ) and $N=2$ (Unit Disc).

The paper begins with some proofs confirming that \eqref{mems_main} exhibits the pull-in instability, \emph{i.e.} there is a number $\eps^{\ast}>0$ such that when $\eps<\eps^{*}$, \eqref{mems_main} has no equilibrium solutions and will touchdown to $u=-1$ in finite time. In \S\ref{sec_numerics}, a moving mesh PDE method (MMPDE) is employed together with an adaptive time stepping scheme to accurately resolve the solution of \eqref{mems_main} very close to touchdown. While touchdown occurs at the origin for certain parameter regimes as in the second order equivalent, it is observed that for $\eps$ below some threshold $\eps_c$, equation \eqref{mems_main} may touchdown at two separate isolated points in the strip case and\jl{, under radially symmetric constraints,} along a ring of points in the unit disc case. Moreover, it is observed that the location of the touchdown set has a dependence on $\eps$ \jl{that can be analyzed}. While multiple touchdown has been observed previously when tailored dielectric profiles $f(x)$ were considered, here the device is uniform ( $f(x)=1$ ) and the location of touchdown can be parameterized through $\eps = \lambda^{-1/2}$. This may potentially allow MEMS devices to perform more exotic tasks or simply extend their lives by spreading wear over a larger area.

In \S\ref{sec_asymptotics}, the location of touchdown for $\eqref{mems_main}$ is analyzed by means of asymptotic expansions which predict that that for the strip case, the two touchdown points are
\begin{equation}\label{intro_tdpoints}
x_c^{\pm} \sim  \pm\Big[1 - \eps^{1/2}f(t_c)^{1/4}[ \eta_0 + f(t_c) \eta_1 + f^2(t_c) \eta_2+\cdots ]\Big], \qquad f(t) = \jl{1 - (1-3t)^{1/3}}
\end{equation}
while for the unit disc \jl{radially symmetric} touchdown occurs on a ring with radius
\begin{equation}\label{intro_tdpoints_disc}
r_c \sim  1 - \eps^{1/2}f(t_c)^{1/4} \eta_{0} -\eps f(t_c)^{1/2} \eta_{\frac14} - \eps^{3/2}f(t_c)^{3/4} \eta_{\frac12}  + \cdots
\end{equation}
where $\eta_0, \eta_1,\eta_2$ and $ \eta_{\frac14}, \eta_{\frac12}$ are numerically determined constants whose values depend on the boundary conditions applied \eqref{mems_main_b}. Note that these asymptotic predictions are in terms of the touchdown time $t_c$ and are valid for $\eps<\eps_c$. In order to estimate the values of $x_c$ and $r_c$, a numerical approximation of $t_c$ is required. These formulae are shown to agree well with full numerics, particularly when $\eps\ll1$. The limiting profile of \eqref{mems_main} as $\inf_{x\in\Omega} u(x,t)\to-1$ is also constructed. In contrast to the quenching profile \eqref{intro_membrane_TD} of the second order problem \eqref{membrane_full}, it is observed that \eqref{mems_main} exhibits a self-similar quenching profile which finalizes to
\begin{equation}\label{intro_bifull_profile}
u(x,t) \to -1 + c_0\left(\frac{|x-x_c|}{\eps^{1/2}}\right)^{4/3}, \qquad \mbox{as} \qquad t\to t_c^{-}
\end{equation}
where the parameter $c_0$ is determined numerically and has value $c_0=0.9060$ for \jl{both} the strip case and touchdown away from the origin in the \jl{radially symmetric} unit disc case. In the unit disc geometry with touchdown at the origin, the numerically obtained value is $c_0 = 0.7265$. The stability of this profile is determined and convergence of the numerical solution of \eqref{mems_main} to the self-similar profile \eqref{intro_bifull_profile} is verified in each case.

\section{Preliminary Results}\label{prelim}

In this section two preliminary results are established. The first result demonstrates that for $\eps$ small enough, \eqref{mems_main} has no equilibrium solution. The second result proves that when no equilibrium solutions exist for \eqref{mems_main}, the solution will touchdown, \emph{i.e.} reach $u(x,t)=-1$, at some point in space in some finite time. These results rely on a positive eigenpair $(\phi_0, \mu_0)$ of the problem
\bsub\label{mems_eigenpair}
\begin{equation}\label{mems_eigenpair_a}
\Delta^2 \phi = \mu \phi, \quad x\in\Omega;
\end{equation}
for the strip and unit disc geometries and the boundary conditions
\begin{equation}\label{mems_eigenpair_b}
\mbox{(Clamped)}: \quad \phi = 0, \quad \partial_n \phi = 0 \quad x\in\partial\Omega; \qquad \mbox{(Navier)}: \quad \phi = 0, \quad \Delta \phi = 0 \quad x\in \partial\Omega
\end{equation}
\esub
In the case of clamped boundary conditions, it is well known that for general two-dimensional geometries, the principal eigenfunction of \eqref{mems_eigenpair} need not be of one sign. Two well known cases are that of the square \cite{CD2} and annulus \cite{CD3}. However, if only the strip and the unit disc are considered, then \eqref{mems_eigenpair} does admit a strictly one signed principal eigenfunction together with a positive eigenvalue. A brief calculation shows that the eigenfunctions for the clamped strip satisfy
\bsub\label{eigenfunction_strip}
\begin{equation}\label{eigenfunction_strip_a}
\phi = C \left[\sin\xi(x-1) -\sinh\xi(x-1) + \left[\frac{\sin2\xi-\sinh2\xi}{\cos2\xi-\cosh2\xi}\right][\cos\xi(x-1) -\cosh\xi(x-1)] \right]\\[5pt]
\end{equation}
where $\xi = \mu^{1/4}$ and 
\begin{equation}\label{eigenfunction_strip_b}
\cos2\xi\cosh2\xi = 1.
\end{equation}
For the clamped unit disc, the eigenfunctions are
\begin{equation}\label{eigenfunction_radial_b}
\phi = C\left[ I_0(\xi r) -\frac{I_0(\xi)}{J_0(\xi) } J_0(\xi r) \right], \qquad J_0(\xi)I_0'(\xi)  = J_0'(\xi)I_0(\xi).
\end{equation}
\esub
where again $\xi = \mu^{1/4}$. In equations \eqref{eigenfunction_strip}, the constant $C$ is fixed by normalization. The case of Navier boundary conditions for this eigenvalue problem were considered in \cite{GUO} where it was shown that $(\phi_{0},\mu_0) = (\phi_{\Omega},\lambda_{\Omega}^2)$ for
\begin{equation}\label{efunctions_lap}
\Delta \phi_{\Omega} + \lambda_{\Omega} \phi_{\Omega} = 0, \quad x\in\Omega; \qquad \phi_{\Omega} = 0. \quad x\in\partial\Omega
\end{equation}
is a positive eigenpair of \eqref{mems_eigenpair_a}. The maximum principle guarantees the positivity of the principal eigenpair of \eqref{efunctions_lap} for any $\Omega\subset\mathbb{R}^N$ \cite{GUO2}. The principal eigenfunction for the strip and unit disc geometries under Navier boundary conditions are therefore,
\bsub\label{eigenfun_Nav}
\begin{align}
\label{eigenfun_Nav_a}\mbox{(Strip)}&: \quad \mu_0 = \frac{\pi^4}{16}  \qquad \phi_0 = C\sin\left(\frac{\pi}{2}(x-1)\right)\\
\label{eigenfun_Nav_b}\mbox{(Unit Disc)}&: \quad \mu_0 = z_0^4 \qquad  \phi_0 = CJ_0(z_0 r),
\end{align}
\esub
where $C$ is a normalization constant and in \eqref{eigenfun_Nav_b} $z_0$ is the first root of $J_0(z_0)=0$. 

The following theorems show that for $\eps$ small enough, equation \eqref{mems_main} admits no equilibrium solutions and will touchdown to $u=-1$ in finite time. The proof techniques involved have been employed previously in \cite{GPW,LW} and rely on a positive eigenfunction of \eqref{mems_eigenpair}. Therefore, in the case of clamped boundary conditions for \eqref{mems_main}, the result is limited to the strip and unit disc geometries.

\noindent {\bf \underline{Theorem 1}}: {\em (c.f. \cite{GPW,LW})  There exists a real $0<\eps^{*}<\infty$ such that for $0<\eps<\eps^{*}$, equation \eqref{mems_main} has no equilibrium solutions when considered on the strip or unit ball  with clamped conditions and any $\Omega\subset\mathbb{R}^2$ for Navier conditions. In addition $\eps^{*} \geq \bar{\eps} =\sqrt{27/4\mu_0}$ where $(\phi_0,\mu_0)$ is a positive eigenpair of \eqref{mems_eigenpair}.}

\noindent {\bf \underline{Proof}}: Take $(\phi_0,\mu_0)$ to be an eigenpair of \eqref{mems_eigenpair} with $\phi_0>0$ and $\mu_0>0$. Multiplying the equilibrium equation of \eqref{mems_main} ( \emph{i.e.} $u_t=0$ ) by $\phi_0$ and integrating gives
\begin{equation}\label{Theorem1_proof1}
\int_{\Omega} \phi_0 \left( \eps^2 \mu_0 u + \frac{1}{(1+u)^2} \right) \, dx = 0
\end{equation}
Clearly \eqref{Theorem1_proof1} cannot hold when the integrand is strictly positive which occurs when the inequality
\begin{equation}\label{Theorem1_proof2}
\eps^2 \mu_0 u + \frac{1}{(1+u)^2} > 0
\end{equation}
is satisfied on $\Omega$. This implies that $\eps^{*}$ is finite. The equality $\eps^2 \mu_0 u = - (1+u)^{-2}$ has exactly one solution when $\eps^2\mu_0 = 27/4$, and no solutions when $\eps^2\mu_0<27/4$.  Therefore, whenever $\eps^2\mu_0<27/4$, \eqref{Theorem1_proof2} holds and \eqref{mems_main} certainly has no equilibrium solutions. Moreover, the smallest positive $\eps$ such that \eqref{mems_main} has an equilibrium solution, $\eps^{\ast}$, satisfies
\begin{equation}\label{Theorem1_proof3}
\eps^{*} \geq \bar{\eps} = \sqrt{\frac{27}{4\mu_0}}.
\end{equation}
Numerical values of $\mu_0$, satisfying the first positive solutions of \eqref{eigenfunction_strip}, $\bar{\eps}$ and $\eps^{\ast}$ are given in Table.~\ref{Table_1} under both clamped and Navier boundary conditions.\endproof

The following theorem shows that \jl{for $\eps < \bar{\eps}$}, when an equilibrium solution to \eqref{mems_main} is not present, touchdown occurs in finite time.

\noindent {\bf \underline{Theorem 2}}: {\em Suppose that $\eps < \bar{\eps} = \sqrt{27/4\mu_0}$, then the solution of \eqref{mems_main} reaches $u = -1$ in some finite time $t_c$ when considered on the strip or unit ball with clamped boundary conditions, or on any bounded $\Omega\subset\mathbb{R}^2$ under Navier boundary conditions. }

\noindent {\bf \underline{Proof}}: The proof follows Theorem 3.1 of \cite{GPW} and relies on the existence of a positive eigenfunction $\phi_0$ of \eqref{mems_eigenpair}. Let $\phi_0$ be normalized by the condition $\int_{\Omega}\phi_0\,dx= 1$. By multiplying \eqref{mems_main_a} by $\phi_0$ and integrating by parts, the equality
\begin{equation}\label{Theorem2_proof1}
\frac{d}{dt} \int_{\Omega} \phi_0 u \, dx = -\eps^2\mu_0 \int_{\Omega} \phi_0 u\,dx - \int_{\Omega} \frac{\phi_0 }{(1+u)^2}\,dx
\end{equation}
is obtained. Defining $E(t) = \int_{\Omega} \phi_0 u\,dx$ where $E(0)=0$ implies that
\begin{equation*}
\frac{dE}{dt} = -\eps^2\mu_0 E - \int_{\Omega} \frac{\phi_0 }{(1+u)^2}\,dx \leq -\eps^2\mu_0 E - \frac{1}{(1+E)^2},
\end{equation*}
where in the last step Jensen's inequality has been applied. Standard comparison principles show that $E(t)\leq F(t)$ where $F(t)$ satisfies
\begin{equation}\label{Theorem2_proof2}
\frac{dF}{dt} = -\eps^2\mu_0 F - \frac{1}{(1+F)^2}, \qquad F(0)=0.
\end{equation}
Equation \eqref{Theorem2_proof2} is separable and so it is solved to show the touchdown time $\bar{t}$, for $F(t)$ at which $F(\bar{t}) = -1$ satisfies
 \begin{equation}\label{Theorem2_proof3}
\bar{t} =  \int_{-1}^{0}\left(\eps^2\mu_0 s + \frac{1}{(1+s)^2}\right)^{-1}\,ds\,.
\end{equation}
The touchdown time for $F(t)$ is finite when this integral converges which occurs when $\eps< \bar{\eps}\equiv \sqrt{27/4\mu_0}$. Finally, since 
 \[
 E(t) = \int_{\Omega} \phi_0 u\, dx \geq  \inf_{x\in\Omega} u \int_{\Omega}\phi_0\, dx =  \inf_{x\in\Omega} u.
 \]
 it follows that 
 \begin{equation}\label{Theorem2_proof4}
  \inf_{x\in\Omega} u \leq E(t) \leq F(t)
 \end{equation}
so that if $\bar{t}$ from \eqref{Theorem2_proof3} is finite, then the touchdown time of \eqref{mems_main}, $t_c$ must also be finite. Therefore when $\eps<\equiv \sqrt{27/4\mu_0}$, $t_c<\bar{t}$ where $\bar{t}$ is given in \eqref{Theorem2_proof3}. In the limit as $\eps\to0^{+}$,  equation \eqref{Theorem2_proof3} has expansion
\begin{equation}\label{Theorem2_proof5}
\bar{t} = \frac{1}{3} + \frac{\eps^2\mu_0}{30} + \mathcal{O}(\eps^4) 
\end{equation}
This provides the asymptotic upper bound on the touchdown time $t_c$ of \eqref{mems_main}
\begin{equation}\label{Theorem2_proof6}
t_c <  \frac{1}{3} + \frac{\eps^2\mu_0}{30}+ \mathcal{O}(\eps^4)
\end{equation}
in the limit as $\eps\to0^{+}$.
\endproof

The preceding analysis demonstrates the presence of the ubiquitous pull-in instability for \eqref{mems_main} for general geometries when the boundary conditions are Navier and for the strip and unit disc domains when clamped boundary conditions are applied. It is an open and challenging problem to prove that \eqref{mems_main} exhibits the pull-in instability for general geometries $\Omega\subset\mathbb{R}^N$ when clamped boundary conditions are applied. A useful byproduct of the analysis presented here is the estimators on the critical pull-in voltage $\eps^{\ast}$ for each of the geometries and boundary conditions considered as collated in Table~\ref{Table_1}.

\begin{table}
\centering
\begin{tabular}{|c|cc|cc|cc|}
\hline
{} &\multicolumn{2}{c|}{$\mu_0$} &\multicolumn{2}{c|}{$\bar{\eps}$}&\multicolumn{2}{c|}{$\eps^{\ast}$}\\
\hline
{} & \mbox{Strip} & \mbox{Unit Disc} &\mbox{Strip} &\mbox{Unit Disc}&\mbox{Strip} &\mbox{Unit Disc} \\[5pt]
\mbox{Navier B.C.s}&$6.0881$ & $33.4452$ &$1.0530$& $0.4492 $&$1.0771$& $0.4695$\\[5pt]
\mbox{Clamped B.C.s}&$31.2852$ & $104.3631$&$0.4645$& $0.2543$&$0.4778$& $0.2683$\\
\hline
\end{tabular}
\vskip0.5cm
\parbox{4.5in}{\caption{Numerical values of the principal eigenvalue $\mu_0$ from \eqref{eigenfunction_strip}, \jl{$\bar\eps$} from Theorem 1,2 and $\eps^{\ast}$ under clamped and Navier boundary conditions. The values of $\eps^{\ast}$ were calculated numerically in \cite{LW} as a saddle-node bifurcation point of equilibrium solutions to \eqref{mems_main}.}\label{Table_1}}
\end{table}

\section{Numerics}\label{sec_numerics}

In order to obtain accurate numerical representations of \eqref{mems_main} close to touchdown, a method which can resolve the rapidly changing spatially localized and temporal features of the equation seems warranted. To facilitate this, the r-adaptive moving mesh scheme MOVCOL4 of \cite{RXW} together with the adaptive time stepping scheme of \cite{BJFW} are implemented. Both schemes take advantage of the underlying invariance of equation \eqref{mems_main} to the transformation 
\begin{equation}\label{numerical_inv}
t\to a t, \quad (1+u) \to a^{1/3}(1+u), \quad x\to a^{1/4}x. 
\end{equation}
A brief overview of the method is now provided, for more details see \cite{BJFW,RXW}. The physical domain is approximated by the grid
\bsub\label{numerical_grid}
\begin{equation}\label{numerical_grid_a}
 x_0 < x_1(t) < \cdots< x_N(t) < x_{N+1},
\end{equation}
the node points of which are evolved with the equation
\begin{equation}\label{numerical_grid_b}
-\gamma X_{t\xi\xi}= (M(X)X_{\xi})_{\xi}.
\end{equation}
\esub
Here $\gamma$ is a small parameter which controls the relaxation timescale to the equidistribution profile, $M(X)$ is known as the monitor function and $ x_i(t)=X(i\Delta \xi,t)$ is a map between the physical domain and a computational domain $\Omega_c = [0,1]$ with coordinate $\xi\in[0,1]$. In calculations, the value $\gamma = 10^{-4}$ was used and the boundary conditions $\dot{X}_0 = \dot{X}_{N+1}=0$ were applied. The monitor function
\begin{equation}\label{numerical_grid_M}
M(X) = \frac{1}{(1+u(X))^3} + \int_{\Omega}\frac{1}{(1+u)^3}\,dx
\end{equation}
was selected which provides a balance between grid points in the region where $M$ is large (\emph{e.g.} where $||1+u||_{\inf}$ is small) and also in regions where the solution is not changing rapidly but modest resolution is still required so that iterative procedures converge. Importantly, with this choice of $M(u)$, equation \eqref{numerical_grid_b} retains the symmetry \eqref{numerical_inv} of the underlying equation. Spatial discretization was effected by a $7^{th}$ order polynomial collocation procedure with evaluation at four Gauss points in each subinterval (c.f. Appendix \ref{AppendixA}). After accounting for boundary conditions, this results in a system of $4(N+1)$ equations for the solution and its first three derivatives at each node point. The mesh equation \eqref{numerical_grid_b} is discretized as follows 
\bsub\label{numerical_CD_M}
\begin{equation}\label{numerical_CD_M_a}
-\gamma \frac{\dot{X}_{i-1} -2\dot{X}_{i}+\dot{X}_{i+1} }{\Delta\xi^2} = \frac{M_{i+\frac12}(X_{i+1}-X_{i}) - M_{i-\frac12}(X_{i}-X_{i-1} )}{\Delta\xi^2}
\end{equation}
where
\begin{equation}\label{numerical_CD_M_b}
M_{i+\frac12} = \frac{M(X_{i+1})+M(X_i)}{2}.
\end{equation}
\esub
The integral term of \eqref{numerical_grid_M} is evaluated by the trapezoid rule on the subintervals defined by the $X_i$s. The efficient simulation of the PDE close to singularity necessitates the use of temporal adaptivity. The underlying symmetry of the problem \eqref{numerical_inv} provides indication on how the time stepping should be adjusted according to the solution magnitude and motivates the introduction of a computational time coordinate \begin{equation}\label{numerical_sulman}
\frac{dt}{d\tau} = g(u)\qquad g(u) = \frac{1}{\inf_{x\in\Omega}||M(u)||}
\end{equation}
where again \eqref{numerical_sulman} retains the underlying symmetry \eqref{numerical_inv} of the underlying problem. The discretized main equation \eqref{mems_main} and equations for the mesh \eqref{numerical_grid_b} are written in terms of the computational time $\tau$ and solved simultaneously as a DAE of form
\begin{equation}\label{numerical_DAE}
 0 = \mathcal{M}(\textbf{y},\tau) \textbf{y}_{\tau} - \textbf{f}(\textbf{y},\tau), \qquad \textbf{y} = ( t(\tau), \textbf{u}, \textbf{X})^{T}.
\end{equation}
Here $\textbf{u}\in\mathbb{R}^{4(N+2)}$ is a vector containing the nodal values of the solution and its first three derivatives while $\textbf{X}\in\mathbb{R}^{N+2}$ is the vector of grid points. The square mass matrix $\mathcal{M}$ is of size $5(N+2)+2$ and has entries filled with the discretizations of \eqref{mems_main} and \eqref{numerical_grid_M} while $\textbf{f}\in\mathbb{R}^{5(N+2)+2}$ represents the discretized right hand sides. The resulting equations are solved in MATLAB with the routine {\tt ode23t}.

\begin{figure}
\centering
\subfigure[$\eps=0.5$]{\includegraphics[width=0.4\textwidth]{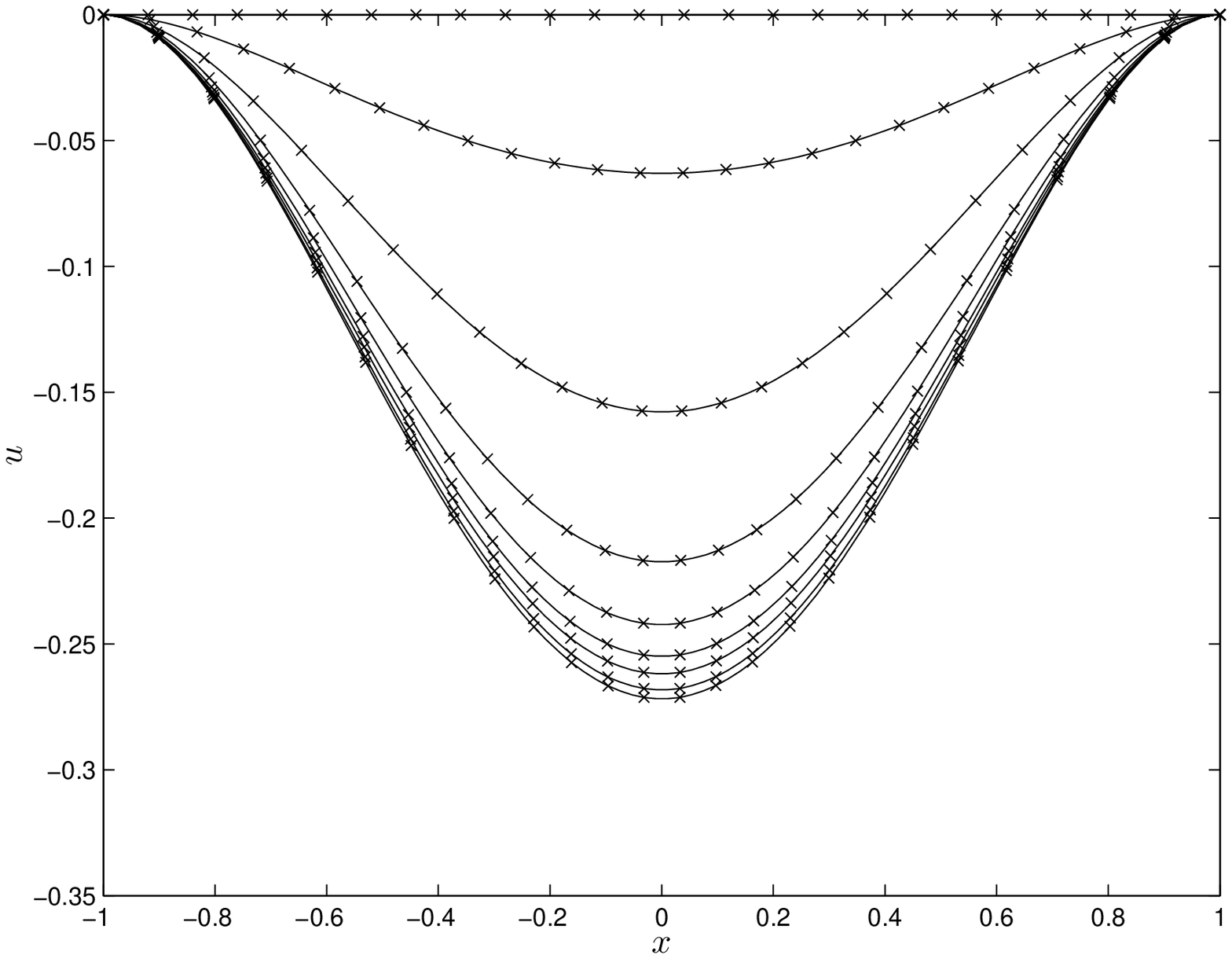}\label{fig_num1a}}\qquad
\subfigure[Touchdown Points]{\includegraphics[width=0.4\textwidth]{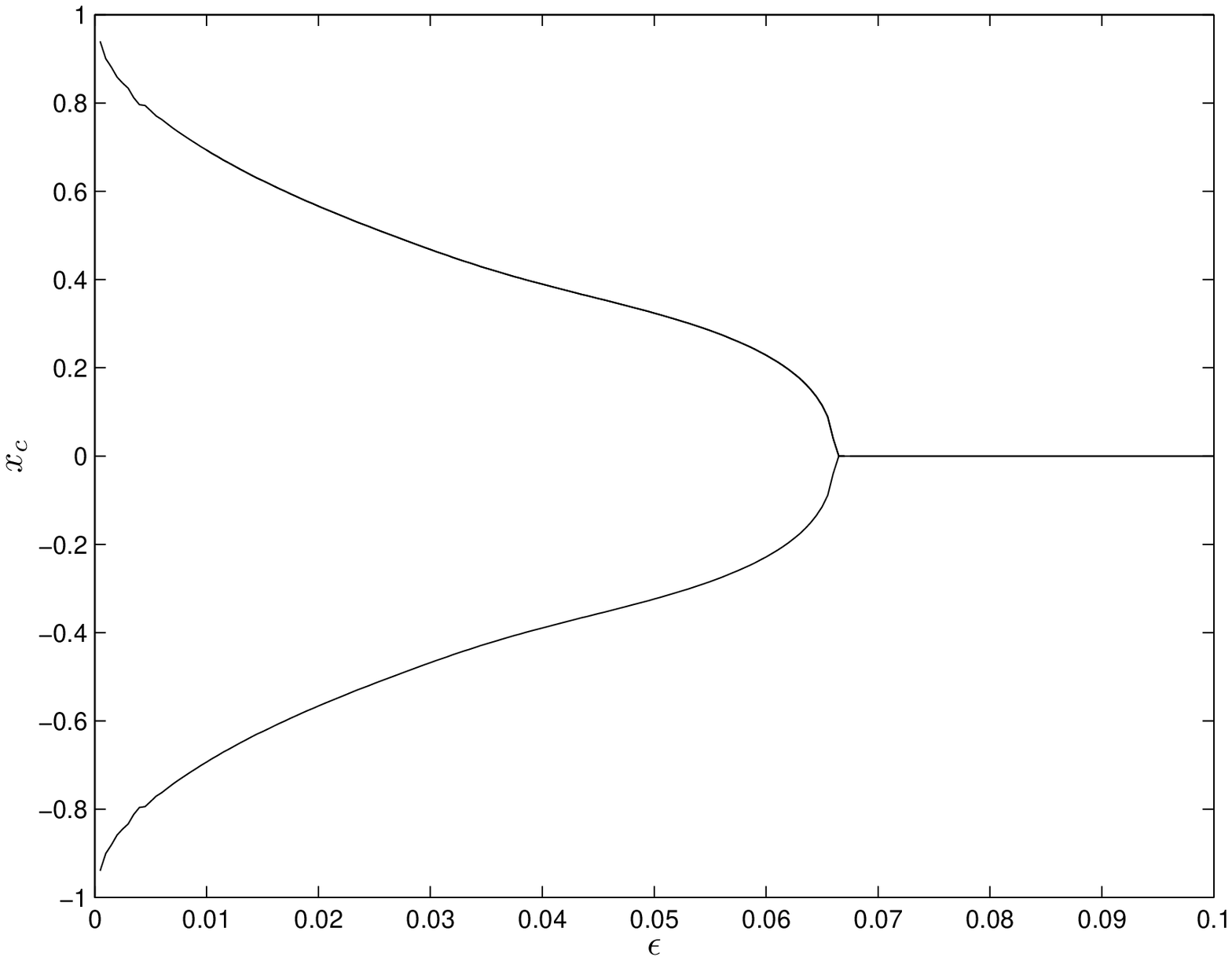}\label{fig_num1d}}
\subfigure[$\eps=0.2$ with $N = 16$]{\includegraphics[width=0.4\textwidth]{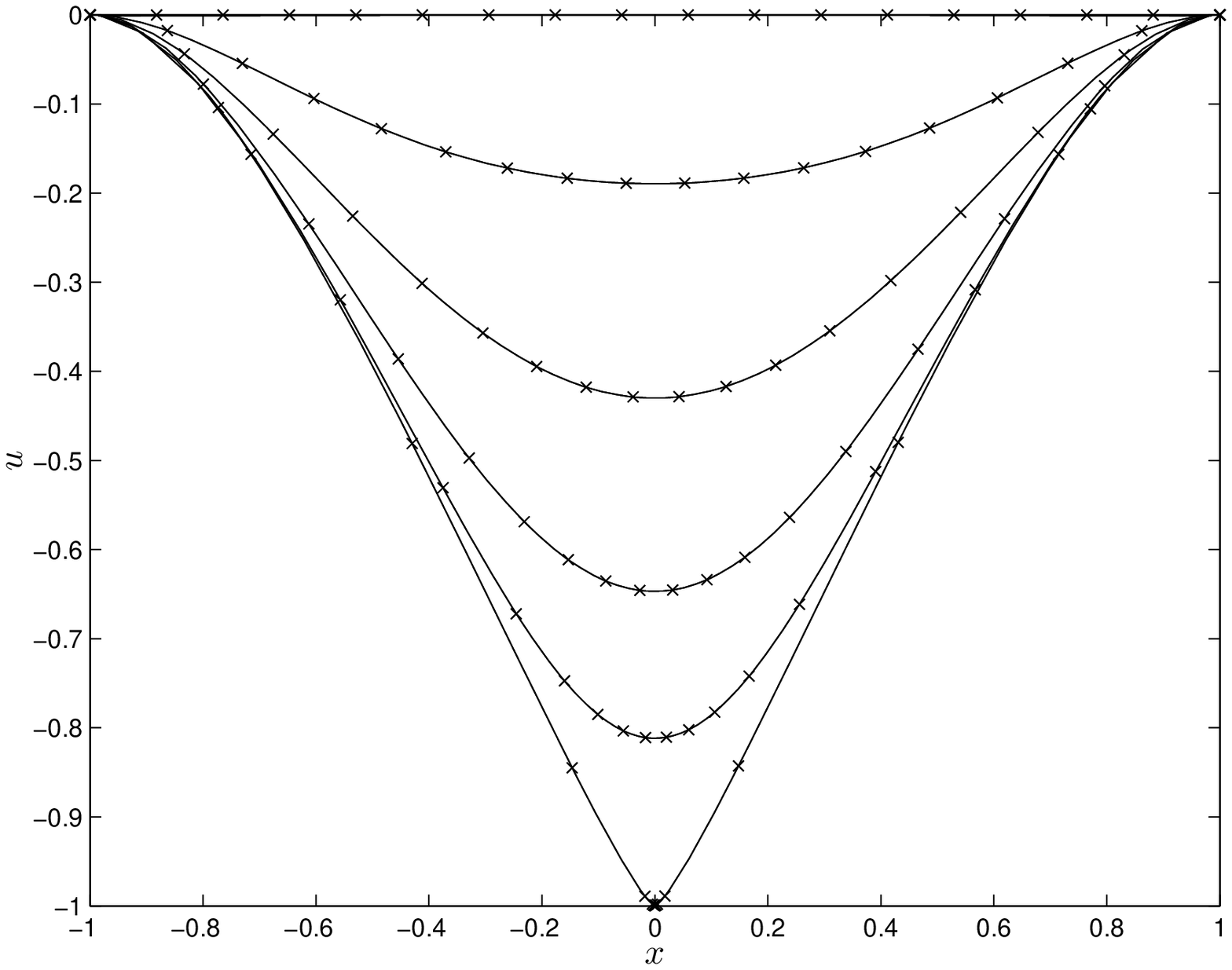}\label{fig_num1b}}\qquad
\subfigure[$\eps=0.2$ with $N = 16$ Zoomed.]{\includegraphics[width=0.4\textwidth]{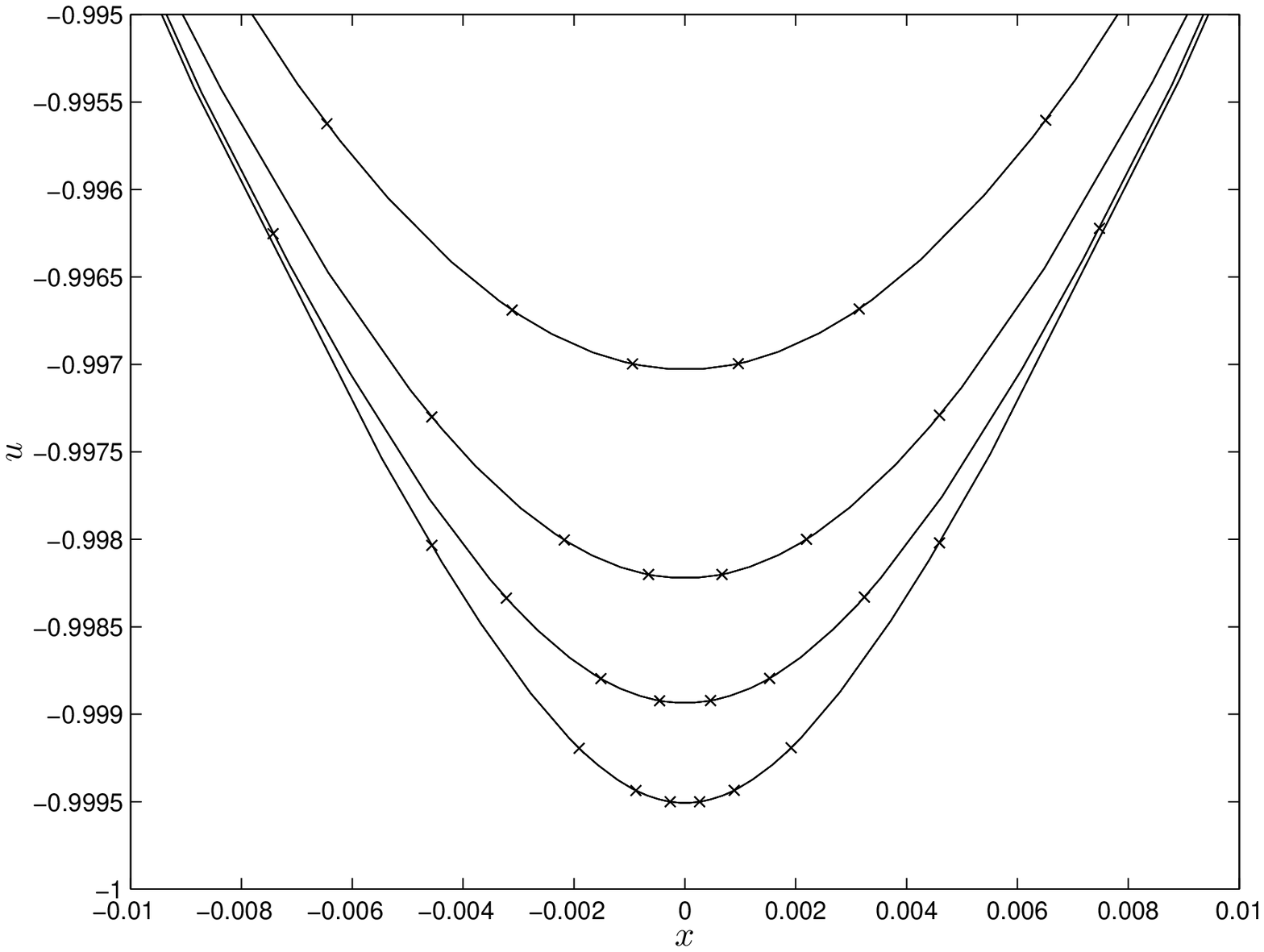}\label{fig_num1b_zoom}}
\subfigure[$\eps=0.02$ and $N=24$.]{\includegraphics[width=0.4\textwidth]{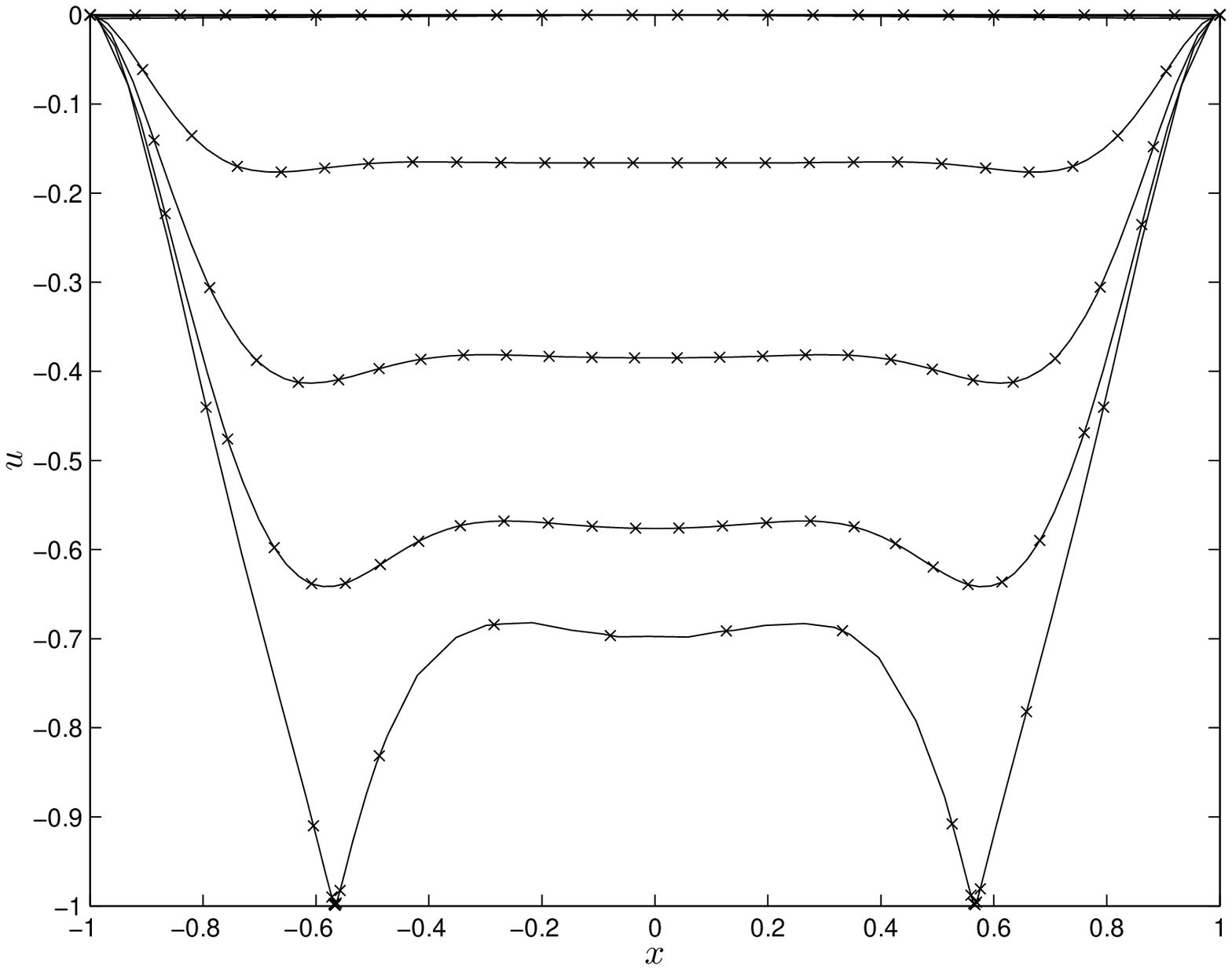}\label{fig_num1c}}\qquad
\subfigure[$\eps=0.02$ and $N=24$ Zoomed.]{\includegraphics[width=0.4\textwidth]{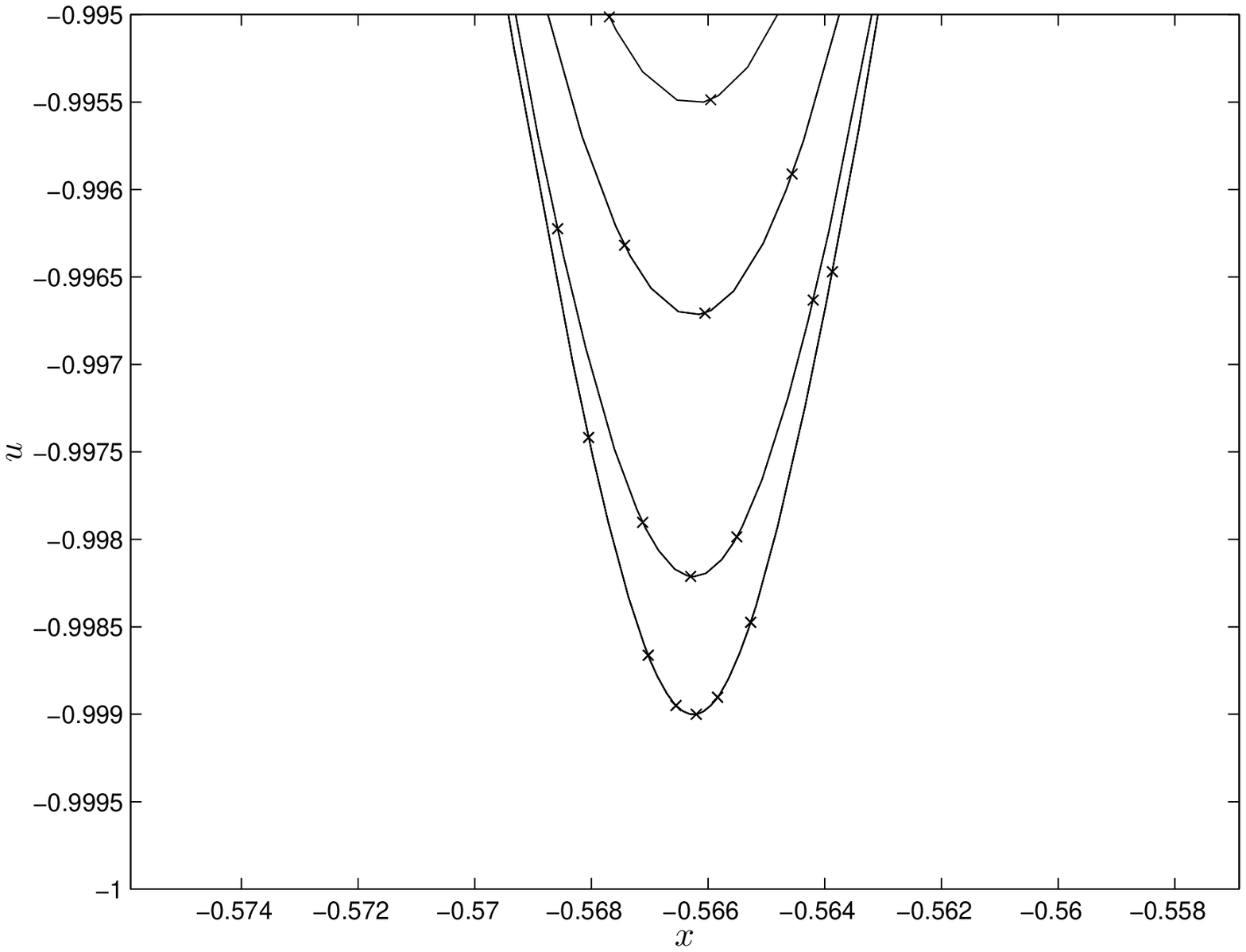}\label{fig_num1c_zoom}}
\parbox{4.5in}{\caption{The above figures relate to numerical solutions of \eqref{mems_main} for the strip domain with clamped boundary conditions. The mesh points are indicated on solutions with small crosses so that their dynamics can be observed. In panel (a), solutions are shown for $\epsilon=0.5>\eps^{*}$ so that touchdown does not occur and a steady state deflection is approached. Panel (b) displays the relationship between touchdown location(s) and the value of $\eps$. The critical value $\eps = \eps_c$, below which touchdown occurs at two points, is approximately $\eps_c = 0.066$. In panels (a) and (c)-(f), solutions are increasing in time from top to bottom. Panel (c) shows solutions for $\epsilon = 0.2<\eps^{*}$ and touchdown is observed at the origin around time $t=0.3833$. In Panel (d) a zoom in of the touchdown region is displayed which shows the refinement of the mesh in this area. In panel (e), solutions are shown for $\eps=0.02$ where touchdown is observed at two separate points, symmetric about the origin around $t=0.3240$. Panel (f) displays a zoom in of the \jl{positive} touchdown region for Panel (e) and again the refinement of the mesh in this region is apparent. \label{fig_num1}}}
\end{figure}

In Fig.~\ref{fig_num1} the three solution regimes for \eqref{mems_main} on the strip under clamped boundary conditions are observed. When $\eps>\eps^{\ast}$ the beam attains a steady equilibrium deflection and does not touchdown (c.f. Fig.~\ref{fig_num1a}). The second solution regime lies in the parameter range $\eps_c<\eps<\eps^{\ast}$ whereby the solution touches down in finite time at the origin only, as displayed in Fig.~\ref{fig_num1c}. \jl{The simulation is halted when $\inf_{x\in\Omega}||1+u(x,t)||$ reaches a specified proximity to $u=-1$. In the case $N=1$ with $\eps = 0.2$, the solution can be followed to $u(0) = -0.99999$ with $t_c-t = \mathcal{O}(10^{-17})$. In the case of multiple touchdown points symmetric about the origin, the solution can be followed to $\inf_{x\in\Omega}||1+u(x,t)|| = -0.999$ where $t_c-t = \mathcal{O}(10^{-10})$. When multiple touchdown points are present it is more challenging to integrate \eqref{mems_main} very close to touchdown as grid points will tend to coalesce on one of the two touchdown points thereby hindering convergence at the other.}ñ On the figures displaying numerical solutions, the grid points are indicated on the curve as crosses and are observed to coalesce on the singularity point as $t\to t_c$ (c.f.  Fig.~\ref{fig_num1b_zoom}). In the third parameter regime $0<\eps<\eps_c$, touchdown occurs in finite time at two isolated points symmetric about the origin (c.f Fig.~ \ref{fig_num1c}) with the location of touchdown as a function of $\eps$ indicated in Fig.~\ref{fig_num1d}. The border of the one and two point touchdown regimes is approximately $\eps_c \approx 0.066$. 

In the radially symmetric unit disc case, touchdown occurs at the origin when $\eps_c<\eps<\eps^{\ast}$ and on an inner ring of points when $\eps<\eps_c\approx 0.075$.

\jl{A possible} interpretation \jl{for this behavior} is that $u=-1$ is an attractor of the system and that the location of touchdown is governed by the critical points of the deflection $u(x,t)$ as the solution enters the basin of attraction for $u=-1$. This would suggest that the source of the multiple touchdown points lies in the dynamics of \eqref{mems_main} for small $t$.

\section{Asymptotics}\label{sec_asymptotics}

\subsection{Small time asymptotics}\label{sec_smalltime}

In this section an analysis of the biharmonic MEMS equation
\begin{equation}\label{mems_main_small}
u_t = -\eps^2 \Delta^2 u - \frac{1}{(1+u)^2}, \quad x\in\Omega; \qquad u(x,0) = 0, \quad x\in\Omega
\end{equation}
is performed in the small time regime $t\to0^{+}$ for strip and disc domains \eqref{mems_main_b} and boundary conditions \eqref{mems_main_c}. In this regime the deflection of the beam is small which allows the $(1+u)^{-2}$ term to be linearized and in this way its influence can be thought of, to leading order, as a uniform forcing term of unit strength.

\begin{figure}
\centering
\includegraphics[width=0.4\textwidth]{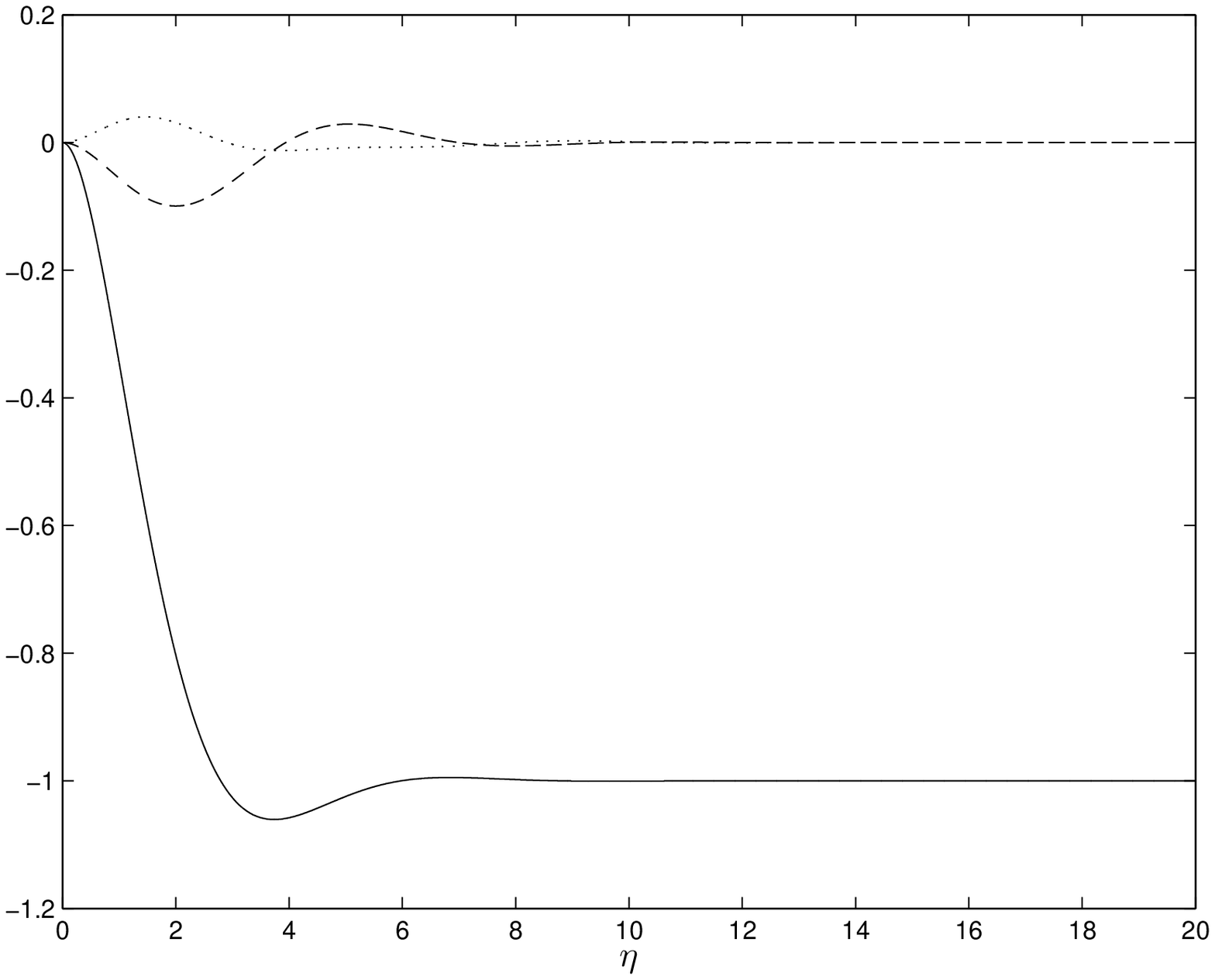}
\includegraphics[width=0.4\textwidth]{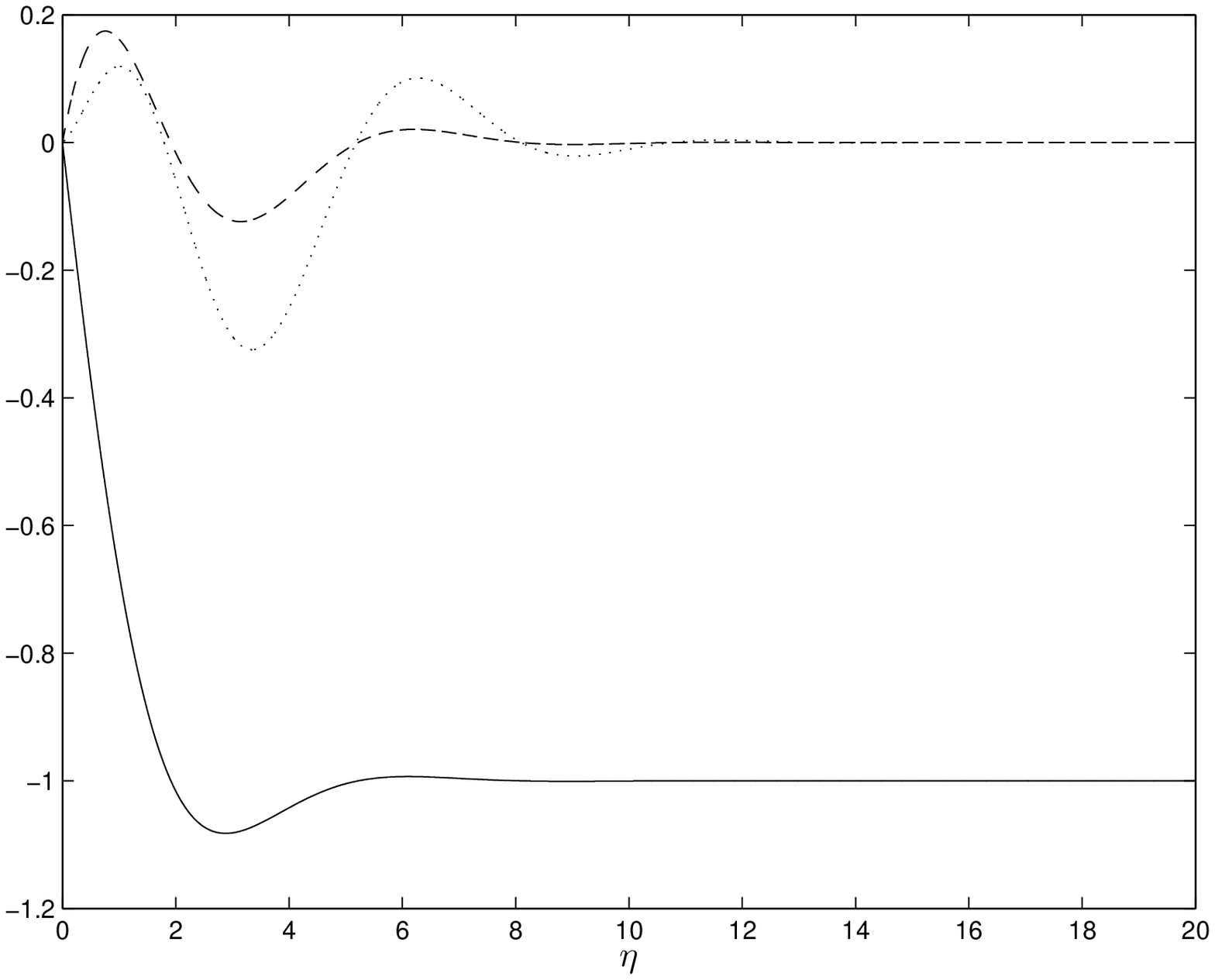}
\parbox{4.5in}{\caption{ Numerical solutions of equations \eqref{mems_full_eqns}. On the left panel $v_0$ (solid curve), $v_1$ (dashed curve) and $v_2$ (dotted curve) are displayed for $N=1$ under clamped boundary conditions. On the right panel,  $v_0$ (solid curve), $v_{\frac14}$ (dashed curve) and $v_{\frac12}$ (dotted curve) are displayed for $N=2$ under Navier boundary conditions.\label{Fig_Profiles} } }
\end{figure}

In a region away from the boundary where the $\eps^2\Delta^2 u$ is negligible for $\eps\ll1$,  the leading order solution satisfies $u(x,t) = \bar{u}(t)$ where
\begin{equation}\label{mems_main_outer}
\bar{u}_t = -\frac{1}{(1+\bar{u})^2}, \quad \bar{u}(0) = 0;  \qquad \bar{u} = -1 + (1-3t)^{1/3}.
\end{equation}
which determines the scale for the solution. This scale together with the scaling invariance \eqref{numerical_inv}, motivates the following expansion for the stretching boundary region in the vicinity of the end point $x=1$
\begin{equation}\label{mems_full_vars}
u(x,t)  = f(t) \, v( \eta,t), \quad \eta = \frac{1-x}{\eps^{1/2}f(t)^{1/4}}, \quad f(t) = 1 - (1-3t)^{1/3}.
\end{equation}
Note that $f(t)= t +\mathcal{O}(t^2)$ as $t\to0$ so an expansion of $v(\eta,t)$ \jl{in powers of $t$} corresponds \jl{at lowest order} to an expansion in small $f(t)$ and matches to the outer region exactly. Employing variables \eqref{mems_full_vars} together with the expansion 
\begin{equation}\label{mems_full_expansion}
v(\eta,t) = \sum_{n=0}^{\infty}f^n(t) v_n(\eta) + \jl{\sum_{k=1, k \ne 4p, p \in \mathbb{N}}^\infty \left(\eps^{1/2}f(t)^{1/4}\right)^k v_{\frac{k}{4}}(\eta)}
\end{equation}
for the solution gives a sequence of problems to be solved for $v_{\frac{k}{4}}(\eta)$, $k=0,1,2,\ldots$. The $\mathcal{O}(\epsilon^{1/2} f(t)^{1/4})$ component of \eqref{mems_full_expansion} is the first correction to the $2r^{-1}u_{rrr}$ term which appears, in the radial case $N=2$, at a lower order due to the expansion not being centred on the origin. \jl{As can be seen from the equations below, when $N=1$, all of the $v_{\frac{k}{4}}$ with non-integer indexes may be chosen zero, so that} the profiles $v_{\jl{\frac{k}{4}}}(\eta)$, for $k\mod4 \neq 0$, play no role in the 1D case. \jl{Equating powers of $f(t)^{1/4}$ yields}
\bsub\label{mems_full_eqns}
\jl{\begin{align}
\label{mems_full_eqns_a}v_{0\eta\eta\eta\eta} -\frac{\eta}{4}v_{0\eta} + v_0 & = -1, \quad  \eta>0;\\[4pt]
v_{1\eta\eta\eta\eta} -\frac{\eta}{4}v_{1\eta} + 2v_1 & = \eps^2(N-1)\, G_1\left(v_0(\eta),v_{\frac{1}{4}}(\eta),v_{\frac12}(\eta),v_{\frac34}(\eta)\right) + \frac{\eta}{2}v_{0\eta}, \quad \eta>0;\label{mems_full_eqns_b} \\[4pt]
v_{2\eta\eta\eta\eta} -\frac{\eta}{4}v_{2\eta} + 3v_2 & = \eps^4(N-1)\, G_2\left(v_0(\eta),v_{\frac14}(\eta),v_{\frac12}(\eta),..., v_{\frac74}(\eta)\right) \nonumber \\[4pt] &\label{mems_full_eqns_c} \ \  -3\left(v_0 - \eta \frac{v_{0\eta}}{4} + v_0^2\right) + \frac{\eta}{2}v_{1\eta}-2v_1, \quad \eta>0;\\[4pt]
\label{mems_full_eqns_d}v_{{\frac14}\eta\eta\eta\eta} -\frac{\eta}{4}v_{{\frac14}\eta} + \frac{5}{4}v_{\frac14} & = 2 (N-1) v_{0\eta\eta\eta}, \quad \eta>0;\\[4pt]
v_{{\frac12}\eta\eta\eta\eta} - \frac{\eta}{4} v_{{\frac12}\eta} + \frac{3}{2} v_{{\frac12}} & = (N-1)\, G_{\frac{1}{2}}\left(v_0(\eta),v_{\frac14}(\eta)\right), \quad \eta>0;\\[4pt]
v_{{\frac34}\eta\eta\eta\eta} - \frac{\eta}{4} v_{{\frac34}\eta} + \frac{7}{4} v_{{\frac34}} & = (N-1)\, G_{\frac34}\left(v_0(\eta),v_{\frac14}(\eta),v_{\frac12}(\eta)\right), \quad \eta>0.
\end{align}}
\jl{In the above, the functions $G_{\frac{k}{4}}$ represent lower order terms that only contribute when $N=2$. In what follows, we retain the first three non-zero terms of expansion (\ref{mems_full_expansion}) when $N=1$, and the first three terms ($v_0$, $v_{\frac14}$ and $v_{\frac12}$) when $N=2$. The above equations are then solved} together with boundary and far field behaviour
\begin{equation}\label{mems_full_eqns_e}
\begin{array}{c}
\mbox{(Clamped)}: \quad v_j(0) = v_{j\eta}(0) = 0, \quad v_{j\eta}, v_{j\eta\eta\eta} \to 0, \quad \eta \to\infty, \quad j= 0,1,2,\jl{\frac{1}{4}},\frac12.\\[15pt]
\mbox{(Navier)}: \begin{array}{c} 
 v_j(0) = v_{j\eta\eta}(0) = 0, \qquad v_{j\eta}, v_{j\eta\eta\eta} \to 0, \quad \eta \to\infty, \qquad j= 0,1,2;\\[5pt]
v_{\frac{k}{4}}(0) = v_{\frac{k}{4}\eta\eta}(0)- v_{\frac{k-1}{4}\eta}(0) = 0, \quad v_{\frac{k}{4}\eta}, v_{\frac{k}{4}\eta\eta\eta} \to 0, \quad \eta \to\infty, \quad k= 1,2
\end{array}
\end{array}
\end{equation}
\esub
The ODEs of \eqref{mems_full_eqns} are solved numerically as boundary value problems on an interval $[0,L]$ with $L$ taken to be sufficiently large so that their limiting behaviour for $\eta \to\infty$ is well manifested. Several profiles $v_0(\eta), v_1(\eta), v_2(\eta)$, $v_{\frac14}(\eta), v_{\frac24}(\eta)$, are displayed in Fig.~\ref{Fig_Profiles} for both boundary conditions.

In the 1D strip case ($N=1$), a solution valid for $x\in(-1,1)$ is obtained by superimposing the left and right boundary phenomena and subtracting the extra far field solution to give a uniform approximation. This gives the small time approximate solution in 1D case
\begin{equation}\label{mems_full_approx_strip}
u(x,t) = f(t)\ds\sum_{n=0}^{2}f^n(t)\left[ v_n\left(\ds\frac{x+1}{\eps^{1/2}f(t)^{1/4}}\right) + v_n\left(\ds\frac{1-x}{\eps^{1/2}f(t)^{1/4}}\right)  \right] - f(t)
\end{equation}

\begin{figure}
\centering
\subfigure[$\eps=0.02$]{\includegraphics[width=0.4\textwidth]{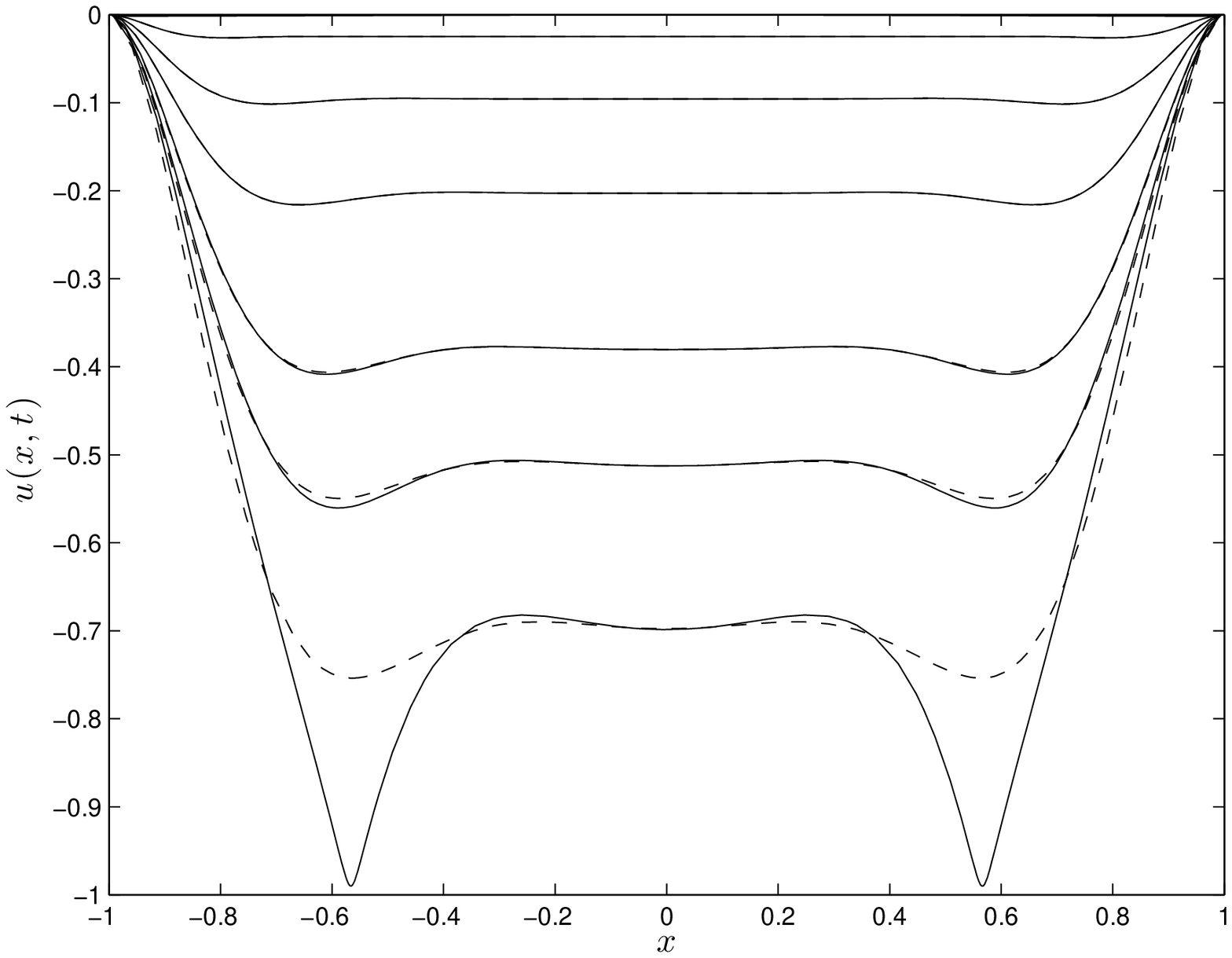}\label{Fig_Numerical_Asy_a}}
\subfigure[$\eps=0.2$]{\includegraphics[width=0.4\textwidth]{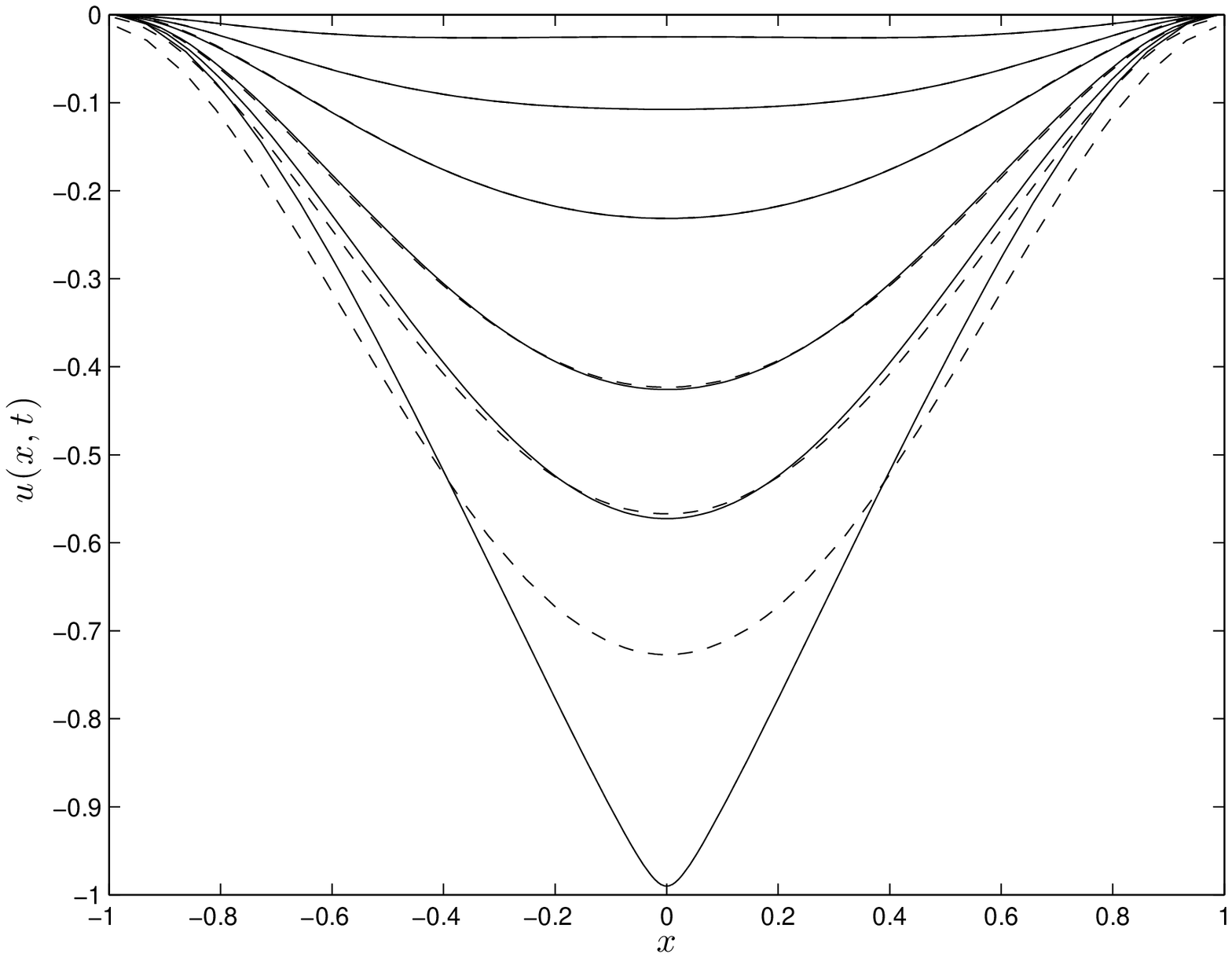}\label{Fig_Numerical_Asy_b}}
\parbox{4.5in}{\caption{ Comparison of full numerical solution ( solid line ) to \eqref{mems_main} on 1D strip with clamped boundary conitions to the asymptotic prediction ( dashed line ) of equation \eqref{mems_full_approx_strip}. Left panel shows $\eps<\eps_c$ so that multiple touchdown points are present while panel (b) has $\eps_c<\eps<\eps^{*}$ so that touchdown occurs at the origin. In both cases, solutions are increasing in time from top to bottom and good agreement between numerics and asymptotics is observed right up till the numerical solution enters the touchdown regime.  \label{Fig_Numerical_Asy} }}
\end{figure}

In Fig.~\ref{Fig_Numerical_Asy}, a comparison of \jl{the} full numerical solution of \eqref{mems_main} and the asymptotic solution \eqref{mems_full_approx_strip} is displayed. Very good agreement is observed for small $t$. As $t\to 1/3^{-}$, $f(t)\to\mathcal{O}(1)$ indicating that the asymptotic solution \eqref{mems_full_approx_strip} breaks down. \jl{Later in time,} a new asymptotic regime based on small $(t_c-t)$ is entered. This touchdown regime is explored in \S\ref{sec_touchdown}.

\begin{figure}
\centering
\subfigure[$\eps=0.02$]{\includegraphics[width=0.4\textwidth]{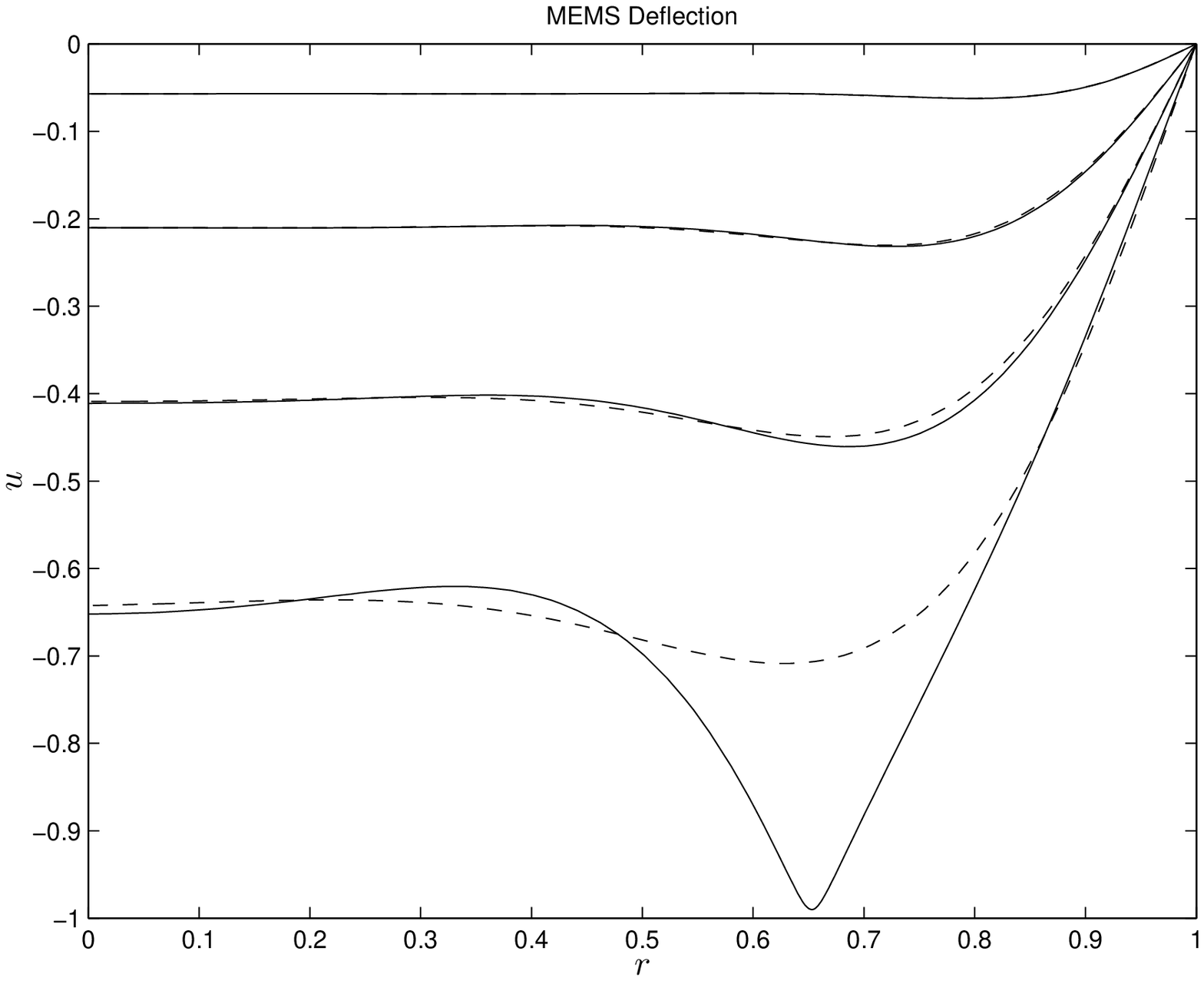}\label{Fig_Radial_NAV_Asy_a}}\qquad
\subfigure[$\eps=0.2$]{\includegraphics[width=0.4\textwidth]{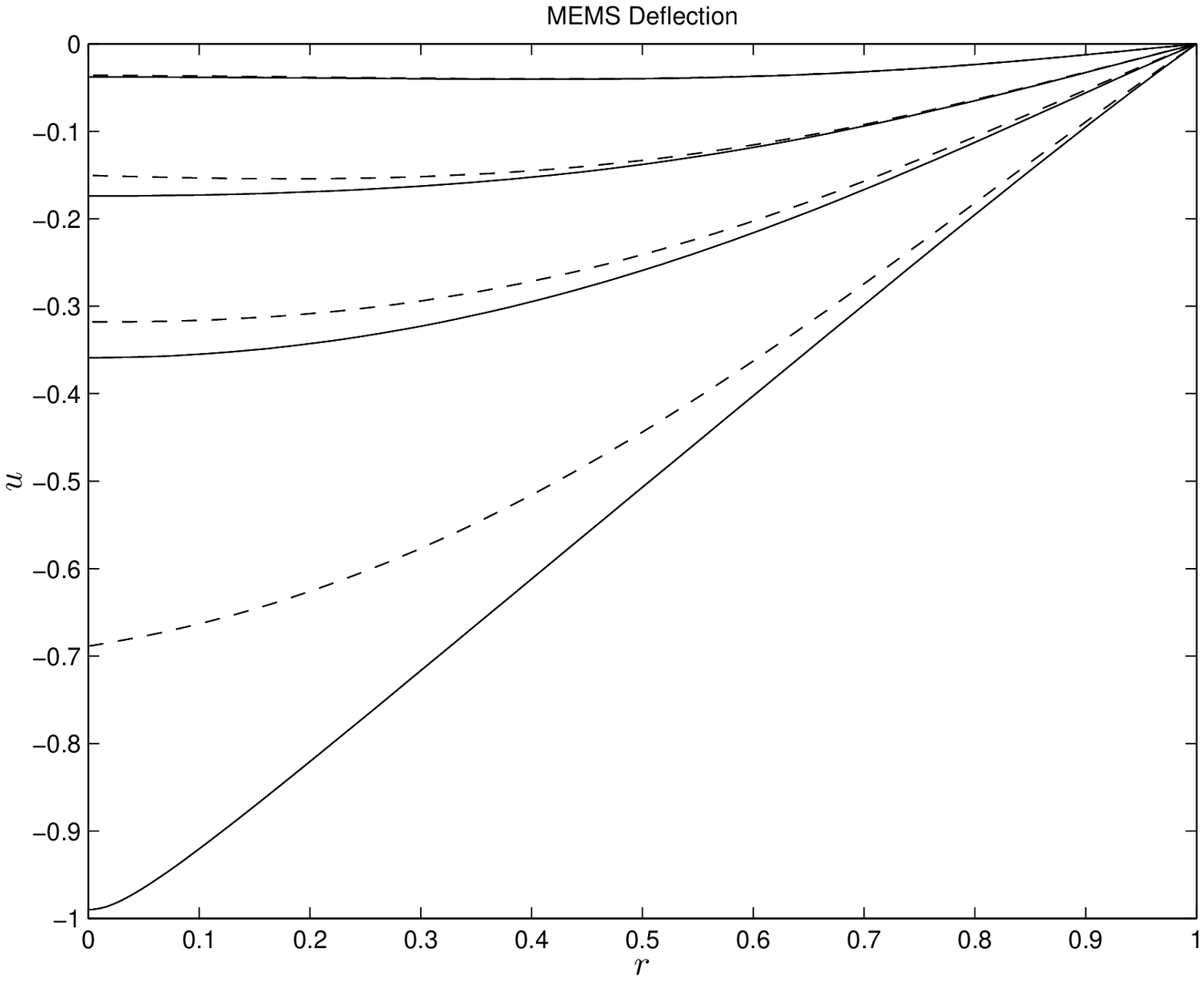}\label{Fig_Radial_NAV_Asy_b}}
\parbox{4.5in}{\caption{ Comparison of full radially symmetric numerical solution (solid line) to \eqref{mems_main} on the unit disc with Navier boundary conditions to the asymptotic prediction (dashed line) of equation \eqref{mems_full_approx_radial}. Left panel shows $\eps<\eps_c$ so that touchdown occurs on a ring of points while panel (b) has $\eps_c<\eps<\eps^{*}$ so that touchdown occurs at the origin. In both cases, solutions are increasing in time from top to bottom and good agreement between numerics and asymptotics is observed right up till the numerical solution enters the touchdown regime.  \label{Fig_Radial_NAV_Asy} }}
\end{figure}

In the unit disk case ($N=2$), the three leading terms in the asymptotic solution, each of which is displayed in Fig.~\ref{Fig_Profiles}  are
\begin{equation}\label{mems_full_approx_radial}
u(r,t) = f(t)\sum_{k=0}^{2} \left(\eps^{1/2}f(t)^{1/4}\right)^k v_{\frac{k}{4}}\left(\ds\frac{1-r}{\eps^{1/2}f(t)^{1/4}}\right).
\end{equation}
Fig.~\ref{Fig_Radial_NAV_Asy} displays a comparison of the full numerical solution to \eqref{mems_main} and the asymptotic solution \eqref{mems_full_approx_radial}. Good agreement is again observed for $t$ small which breaks down as $t\to1/3$ and the touchdown regime is entered.

\subsubsection{Estimation of Touchdown Points}

To estimate the touchdown points of \eqref{mems_main}, the critical points of the small $t$ approximations \eqref{mems_full_approx_strip} and \eqref{mems_full_approx_radial} are examined. The first trough of the profile $v(\eta,t)$, defined in \eqref{mems_full_vars}, serves as an estimator of the touchdown points and so its approximate value is determined asymptotically from \eqref{mems_full_expansion}. \jl{Minima of $u$} are candidates for touchdown points with their location determined by the zeros of the derivative of $u$ where, to leading order, the candidates for the touchdown points satisfy

\bsub\label{Touchdowncriticalzero}
\begin{align}
\label{Touchdowncriticalzero_a} v_0'\left(\ds\frac{1+x_c}{\eps^{1/2}f(t)^{1/4}}\right)&=v_0'\left(\ds\frac{1-x_c}{\eps^{1/2}f(t)^{1/4}}\right), &(N=1)\\[4pt]
\label{Touchdowncriticalzero_b} v_0'\left(\ds\frac{1-r_c}{\eps^{1/2}f(t)^{1/4}}\right) &= 0 & (N=2)
\end{align}
\esub

Note that when $1+x_c = \mathcal{O}(\eps^{1/2}f(t)^{1/4})$, $1-x_c = \mathcal{O} (1)$ and so for $\eps^{1/2}f(t)^{1/4}\ll1$, a zero of the left hand side of \eqref{Touchdowncriticalzero_a} corresponds to the far field, \emph{i.e.} flat, region of the right hand side of \eqref{Touchdowncriticalzero_a}. In other words, for $\eps^{1/2}f(t)^{1/4}\ll1$ in the strip case $N=1$, the two propagating regions do not interact directly and the critical points are the local maximums of the profile $v(\eta,t)$ inside each of the two regions. This assumption breaks down when $\eps^{1/2}f(t)^{1/4}=\mathcal{O}(1)$ as the two waves will superimpose to generate more complex solutions of \eqref{Touchdowncriticalzero}.

\begin{figure}
\centering
\includegraphics[width=0.4\textwidth,clip]{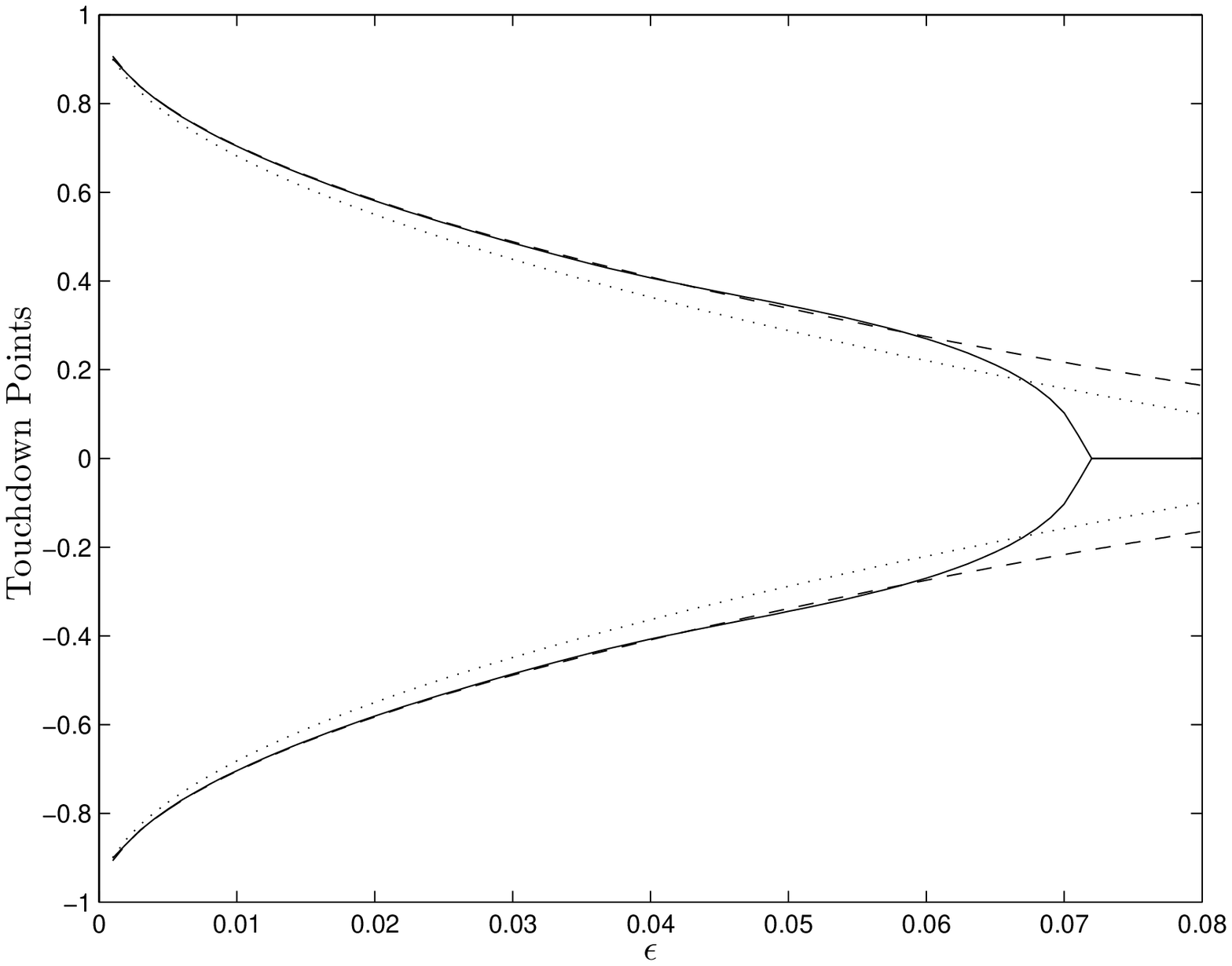}\qquad
\includegraphics[width=0.4\textwidth,clip]{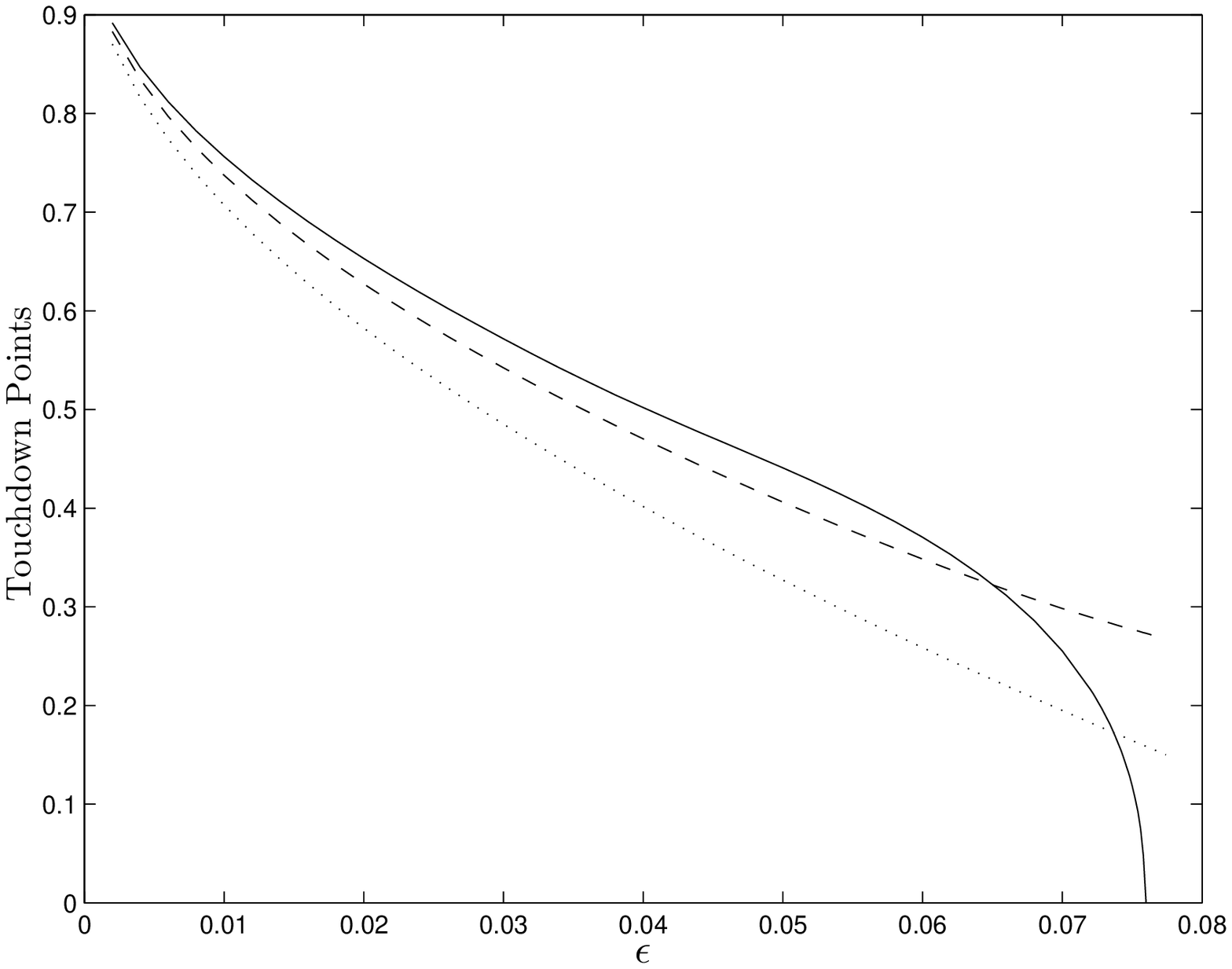}
\parbox{4.5in}{\caption{Touchdown location for \eqref{mems_main} from full numerics ( solid line ), compared with asymptotic formula \eqref{TD_asy} with $t_c$ from full numerics ( dashed line ) and asymptotic formula \eqref{TD_asy} with $t_c=1/3$ ( dotted line ). Left figure case $N=1$ with clamped boundary conditions, right figure case $N=2$ with Navier boundary conditions. \label{Fig_Touchown_Predict} }}
\end{figure}

The critical point inside each expanding region, $\eta_c(t)$, satisfies
\[
\begin{array}{rl} (N =1) &\eta_c(t) = \eta_0 + f(t) \eta_1 + f^2(t) \eta_2 + \cdots\\[5pt]
(N=2) &\eta_c(t) = \eta_0 + \epsilon^{1/2}f(t)^{1/4} \eta_{\frac14} + \epsilon f(t)^{1/2}\eta_{\frac12} + \cdots\end{array} \qquad f(t) = 1 - (1-3t)^{1/3},
\]
where the corrections are determined asymptotically from the condition $v_{\eta}(\eta_c(t),t)=0$. In the $N=1$ case, this provides the condition
\begin{align*}
0 &= v_{0\eta}(\eta_c) + f v_{1\eta}(\eta_c) + f^2 v_{2\eta}(\eta_c) + \cdots\\[5pt]
{} &= v_{0\eta}(\eta_0) + f [ \eta_1 v_{0\eta\eta}(\eta_0) + v_{1\eta}(\eta_0) ] \\[5pt]
{} &+ f^2 [ v_{2\eta}(\eta_0) +\eta_2v_{0\eta\eta}(\eta_0) + \eta_1v_{1\eta\eta}(\eta_0) + \frac{\eta_1^2}{2} v_{0\eta\eta\eta}(\eta_0) ] + \cdots 
\end{align*}
which gives the following definition for the corrections $\eta_j$, $j= 0,1,2;$
\begin{equation}
\begin{array}{c}
v_{0\eta}(\eta_0) = 0, \qquad \eta_1 = -\ds\frac{v_{1\eta}(\eta_0)}{v_{0\eta\eta}(\eta_0)}, \\[10pt] \eta_2 = \ds\frac{-1}{v_{0\eta\eta}(\eta_0)}\left[   v_{2\eta}(\eta_0) + \eta_1v_{1\eta\eta}(\eta_0) + \ds\frac{\eta_1^2}{2} v_{0\eta\eta\eta}(\eta_0)\right]. \end{array}
\end{equation}
A similar calculation can be performed for the $N=2$ case and so the values of $\eta_0,\eta_1, \eta_2$ for $N=1$ with clamped boundary conditions and $\eta_0,\eta_{\frac14}, \eta_{\frac12}$ for $N=2$ with Navier boundary conditions are found to be
\bsub\label{onstants_xc}
\begin{align}
\label{constants_xc_clamped} \mbox{(Clamped)}: \qquad &\eta_0 = 3.7384 , \qquad \eta_1 = -0.6641, &\eta_2 = 0.1085\\[5pt]
\label{constants_xc_Navier} \mbox{(Navier)}: \qquad &\eta_0 = 2.8832 , \qquad \eta_{\frac14}= \phantom{-}0.3533, &\eta_{\frac12} = 0.9457.
\end{align}
\esub
This now allows for the two critical points $x_c^{\pm}(t)$ in the strip case $N=1$ and the ring of touchdown points $r_c(t)$ in the \jl{radially symmetric} unit disc case $N=2$ to be specified as
\bsub\label{TD_asy}
\begin{align}
\label{TD_asy_N=1} N=1& \quad x_c^{\pm}(t) =  \pm\Big[1 - \eps^{1/2}f(t)^{1/4}[ \eta_0 + f(t) \eta_1 + f^2(t) \eta_2 ]\Big], \\
\label{TD_asy_N=2} N=2& \quad r_c  =  1 - \eps^{1/2}f(t_c)^{1/4} \eta_{0} -\eps f(t_c)^{1/2} \eta_{\frac14} - \eps^{3/2}f(t_c)^{3/4} \eta_{\frac12}  + \cdots
\end{align}
\esub
Note that the approximation for the touchdown locations requires $t_c$, the touchdown time of \eqref{mems_main}. As observed in Fig.~\ref{Fig_Touchown_Predict}, asymptotic formula \eqref{TD_asy} captures the location of touchdown very well, particularly when $\eps\ll\eps_c$. As $\eps\to\eps_c^{-}$, the approximation breaks as the left and right boundary effects are superimposing and so the touchdown points are no longer simply the minima of the isolated profile $v(y,t)$.ali

\section{Touchdown regime}\label{sec_touchdown}

To establish a blow up profile in the touchdown regime, the techniques of \cite{BW1} are employed. The correct similarity variables are investigated by initially rescaling equation \eqref{mems_main_a} with 
\[
u = -1 + U \hat{u}(\hat{x},\hat{t}), \qquad t = T \hat{t}, \quad x = L\hat{x}
\]
which results in
\[
\frac{U}{T} \hat{u}_{ \hat{t}} = -\eps^2 \frac{U}{L^4} \Delta^2_{\hat{x}}\hat{u} - \frac{1}{U^2 \hat{u}^2}.
\]
A balance of all terms suggests scaling with $L\sim T^{1/4}$ and $U\sim{T}^{1/3}$ and so an appropriate self similar solution would be of form
\[
u = -1 + R(t)^{1/3} v\left(\frac{x}{R(t)^{1/4}}\right)
\]
where $R(t)$ is the quenching rate of the solution. In general, rigorous determination of $R(t)$ is a difficult problem and so we make reasonable guesses and investigate their validity with numerical calculations. This approach is not definitive, however, as the case of blow-up in critical NLS\cite{Gadi}, whereby the rate has been found to satisfy the so called \emph{loglog} law
\[
R(t) \sim  \frac{2\pi(t_c-t)}{\log(-\log(t_c-t))}, \qquad t\to t^{-}_c,
\]
indicates. Numerical verification of this rate law would require accurate solutions for almost surely unobtainably small values of $(t_c-t)$. As such, the evidence presented here for self-similar quenching awaits rigorous verification. The case of quenching solutions in the strip and unit disc geometries are treated separately and both appear to be self-similar in nature.

\subsection{ Touchdown solutions in 1D}

In similarity variables
\begin{equation}\label{td_var}
u(x,t) = -1 + (t_c-t)^{1/3}v(\eta,s), \qquad \eta = \frac{x-x_c}{\eps^{1/2}(t_c-t)^{1/4}}, \quad s= -\log(t_c-t)
\end{equation}
\eqref{mems_main_a} is transformed to
\begin{equation}\label{main_v}
v_s = -v_{\eta\eta\eta\eta} -\frac{\eta}{4}v_{\eta} + \frac{v}{3} - \frac{1}{v^2}; \qquad (\eta,s) \in \mathbb{R}\times\mathbb{R}^{+}.
\end{equation}
Far field and initial conditions for $v(\eta,s)$ are now discussed. The behaviour of $v(\eta,s)$ for $\eta \to\pm\infty$ corresponds to a solution of $u(x,t)$ for $x\neq x_c$ as $t\to t_c^{-}$. Assuming a localized quenching solution at $x=x_c$, it can be expected that $u_t = \mathcal{O}(1)$ in a region away from $x_c$ as $t\to t_c^{-}$. Now,
\begin{equation}\label{main_v_bc}
u_t = (t_c -t)^{-2/3} \left[ v_s + \frac{\eta}{4} v_{\eta} -\frac{v}{3}\right]
\end{equation}
and so the condition that $u_t = \mathcal{O}(1)$ implies that 
\begin{equation}\label{main_v_bc2}
v_s + \frac{\eta}{4} v_{\eta} -\frac{v}{3} = \mathcal{O}((t_c -t)^{2/3}), \qquad t \to t_c^{-}
\end{equation}
For a fixed $x\neq x_c$, the limit $t \to t_c^{-}$ corresponds to $|\eta|\to\infty$ and so \eqref{main_v_bc2} augments \eqref{main_v} to establish
\bsub\label{main_v2}
\begin{align}
\label{main_v2_a} v_s &= -v_{\eta\eta\eta\eta} -\frac{\eta}{4}v_{\eta} + \frac{v}{3} - \frac{1}{v^2}, \qquad (\eta,s) \in \mathbb{R}\times\mathbb{R}^{+};\\[5pt]
\label{main_v2_b} v_s &= \frac{v}{3} -\frac{\eta}{4} v_{\eta}, \qquad \eta \to\pm \infty;
\end{align}
\esub

A key step is to determine the limiting behavior of solutions to \eqref{main_v2} for any fixed $\eta$ as $s\to\infty$. One obvious candidate for an equilibrium state is the constant $\bar{v} = 3^{1/3}$. An analysis of its stability leads one to consider the eigenvalue problem
\begin{equation}\label{membrane_mems4}
\mathcal{L}_2 w = \jl{\mu} w, \qquad \mathcal{L}_m \equiv -\left(-\frac{d^2}{d\eta^2}\right)^{m} - \frac{\eta}{4}\frac{d}{d\eta} + I.
\end{equation} 
for $m=2$. The spectrum of the operator $\mathcal{L}_m$ in the weighted space $L^2_{\rho}(\mathbb{R})$ where $\rho= e^{-a|\eta|^{\nu}}$, with $a$ positive  ( c.f. \cite{BGW,GJFW}), is
\begin{equation}\label{membrane_mems5}
\sigma(\mathcal{L}_m) = \left\{ \jl{\mu}_k = 1-\frac{k}{m}; \quad k = 0,1,2,\ldots \right\}
\end{equation}
and so there are two linearly unstable modes associated with this equilibrium, $\jl{\mu}_0 = 1$ and $\jl{\mu}_1 = 1- 1/m$ for $m\geq2$. The instability associated with the mode $\jl{\mu}_0 = 1$ is generated by the invariance of the touchdown time $t_c$ and is therefore not a true instability. The instability associated with the $\jl{\mu}_1 = 1-1/m$ mode represents a true instability when $m\geq 2$. The presence of this positive eigenvalue indicates that equation \eqref{main_v2} does not satisfy $v(\eta,s)\to 3^{1/3}$ for fixed $\eta$ as $s \to \infty$, thus we seek equilibrium solutions  of the following nonlinear problem
\bsub\label{main_v3}
\begin{align}
\label{main_v3_a}  \bar{v}_{\eta\eta\eta\eta} +\frac{\eta}{4}\bar{v}_{\eta} - \frac{\bar{v}}{3} + \frac{1}{\bar{v}^2}&=0, \qquad -\infty<\eta<\infty;\\[5pt]
\label{main_v3_b}  \frac{\bar{v}}{3} -\frac{\eta}{4} \bar{v}_{\eta}&=0, \qquad \eta \to\pm\infty.
\end{align}
\esub
 and investigate the multiplicity and stability of its solution. 
 
 The robin condition of \eqref{main_v3_b}  suggests that \eqref{main_v3} admits a far field series solution of form
\begin{equation}\label{eq_v_1}
\bar{v}(\eta) \sim v_p \equiv \sum_{n=0}^{\infty} c_n|\eta|^{4/3 -4n} \qquad |\eta| \to\pm\infty.
\end{equation}
Here the constants $c_n=c_n(c_0)$ are functions of the parameter $c_0$ for $n\geq1$ and can be determined by lengthy but straightforward manipulations, \emph{e.g.} $c_1 = 40c_0/81 + c_0^{-2}$. The parameter $c_0$ plays the role of a nonlinear eigenvalue and it is expected that \eqref{main_v3} will have solutions for isolated values only. Taking the limit $t\to t_c^{-}$ for fixed $x\neq x_c$ corresponds to the limit $|\eta|\to\infty$ and therefore, in physical co-ordinates, the touchdown profile is expected to satisfy
\begin{equation}\label{touchdown_profile}
u(x,t) \sim -1 + c_0\left|\frac{x-x_c}{\eps^{1/2}}\right|^{4/3} + c_1 (t_c-t)\left|\frac{x-x_c}{\eps^{1/2}}\right|^{-8/3} +\cdots \qquad \mbox{as} \qquad t\to t_c^{-}
\end{equation}
Additional boundary conditions are now obtained for \eqref{main_v3} by suppressing exponentially growing modes of the linearization of \eqref{main_v3} about $\bar{v}$ for large $\eta$. To analyze linearized perturbations of \eqref{main_v3} about $v_p$, set $\bar{v} = v_p + \sigma w$ where $\sigma\ll1$ to arrive at the equation
\begin{equation}\label{main_eq_w}  w_{\eta\eta\eta\eta} +\frac{\eta}{4}w_{\eta} - \frac{w}{3} - 2 \frac{w}{v_p^3}=0, \qquad -\infty<\eta<\infty.
\end{equation}
For large $|\eta|$, a WKB anzatz solution of form 
\[
w\sim \exp\left[ \frac{1}{\delta} \sum_{k=0}^{\infty} \delta^k g_k(\zeta) \right] , \qquad  \eta = \frac{\zeta}{\nu},
\]
for $\nu\ll1$ and $\delta = \nu^{4/3}$  produces the leading order equation
\[
g_{0\zeta}^4 + \frac{\zeta}{4} g_{0\zeta} = 0
\]
which admits three exponential solutions
\[
g_{0j}(\zeta) = -3|\zeta|^{4/3} 2^{-8/3}  \exp\left[{\frac{2\pi i j}{3}}\right], \quad j=0,1,2.
\]
The terms $\exp(g_{0j})$ for $j=1,2$ are growing as $\eta\to\pm\infty$ and need to be suppressed in the solution of \eqref{main_v3}. The mode corresponding to $g_0'=0$ is $w= \eta^{4/3}$ which represents an arbitrary change in the value of $c_0$. At the following order $\exp(g_1(\zeta)) = \eta^{-10/9}$ which now gives the following full specification for $\bar{v}(\eta)$
\bsub\label{main_v4}
\begin{align}
\label{main_v4_a}  \bar{v}_{\eta\eta\eta\eta} +\frac{\eta}{4}\bar{v}_{\eta} - \frac{\bar{v}}{3} + \frac{1}{\bar{v}^2}&=0, \quad -\infty<\eta<\infty.\\[5pt]
\label{main_v4_b}  \bar{v}\sim \sum_{n=0}^{\infty} c_n\eta^{4/3 -4n} &+ \bar{C}|\eta|^{-10/9}\exp\left[-3|\eta|^{4/3} 2^{-8/3} \right], \qquad \eta \to \pm\infty
\end{align}
\esub

Extracting information from \eqref{main_v4} is analytically challenging as it involves solving a fourth order, nonlinear, non constant co\-efficient and non-variational differential equation. This motivates the use of numerical techniques to analyze the multiplicity and stability of solutions to \eqref{main_v4}. 


\subsubsection{Numerical and stability analysis}

This section deals with the numerical determination and linear stability of solutions to \eqref{main_v4}. Related similarity ODES have been solved by several authors in the context of pinch off dynamics for thin films \cite{BW2,BW3} and a framework for their solution is well established. Equation \eqref{main_v4_a} is solved by first applying a centered difference discretization scheme to the derivative terms on a uniform grid of $[-L,L]$. The Robin condition \eqref{main_v3_b} is discretized and applied to remove the ghost points from both end points and thus effectively yields four boundary conditions for the system. The application of the Robin condition at two nodal points enforces the far field behaviour $\bar{v}\sim c_0|\eta|^{4/3}$ and also eliminates exponentially growing terms.

This discretization leads to a large system of nonlinear equations to be solved via a relaxed Newton's Method \cite{AMR}. The iterations are initialized with a solution of the reduced equation 
\begin{equation}\label{num_init}
\frac{\eta}{4}\tilde{v}_{\eta} - \frac{\tilde{v}}{3} + \frac{1}{\tilde{v}^2}=0  \qquad \tilde{v} = \sqrt[3]{c_0^3\eta^4+ 3}, \qquad c_0>0
\end{equation}
 over a wide range of positive $c_0$ until convergence is achieved. This initial guess has the advantage of satisfying the far field behaviour exactly for a given $c_0$ and also being smooth at the origin. The size \jl{$L$ of the system} is taken to be sufficiently large so that the far field behaviour is well manifested.
\begin{figure}
\centering
\subfigure[Profile $\bar{v}_1(\eta)$]{\includegraphics[width=0.4\textwidth]{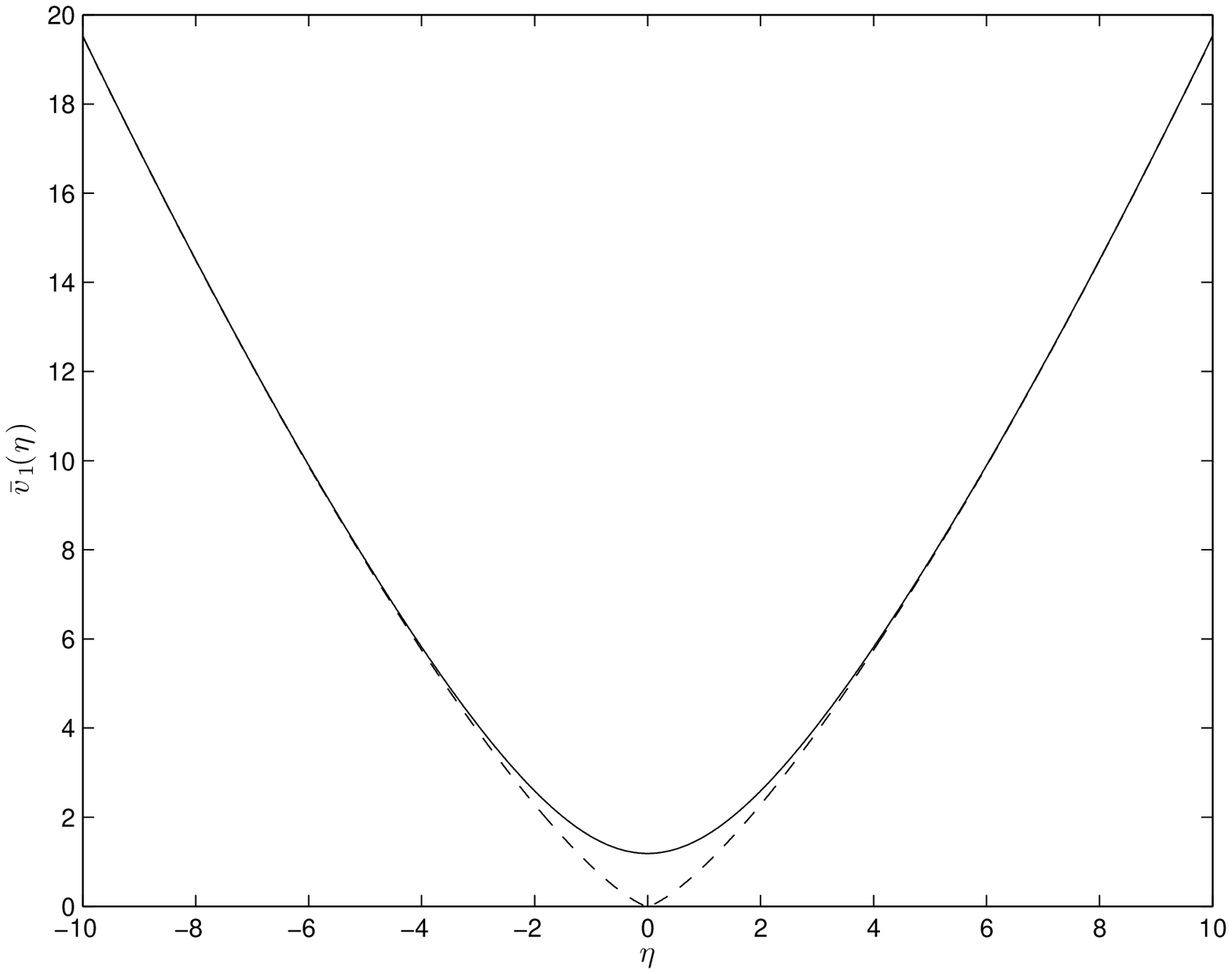}\label{fig_num_TDa}}\qquad
\subfigure[Profile $\bar{v}_2(\eta)$]{\includegraphics[width=0.4\textwidth]{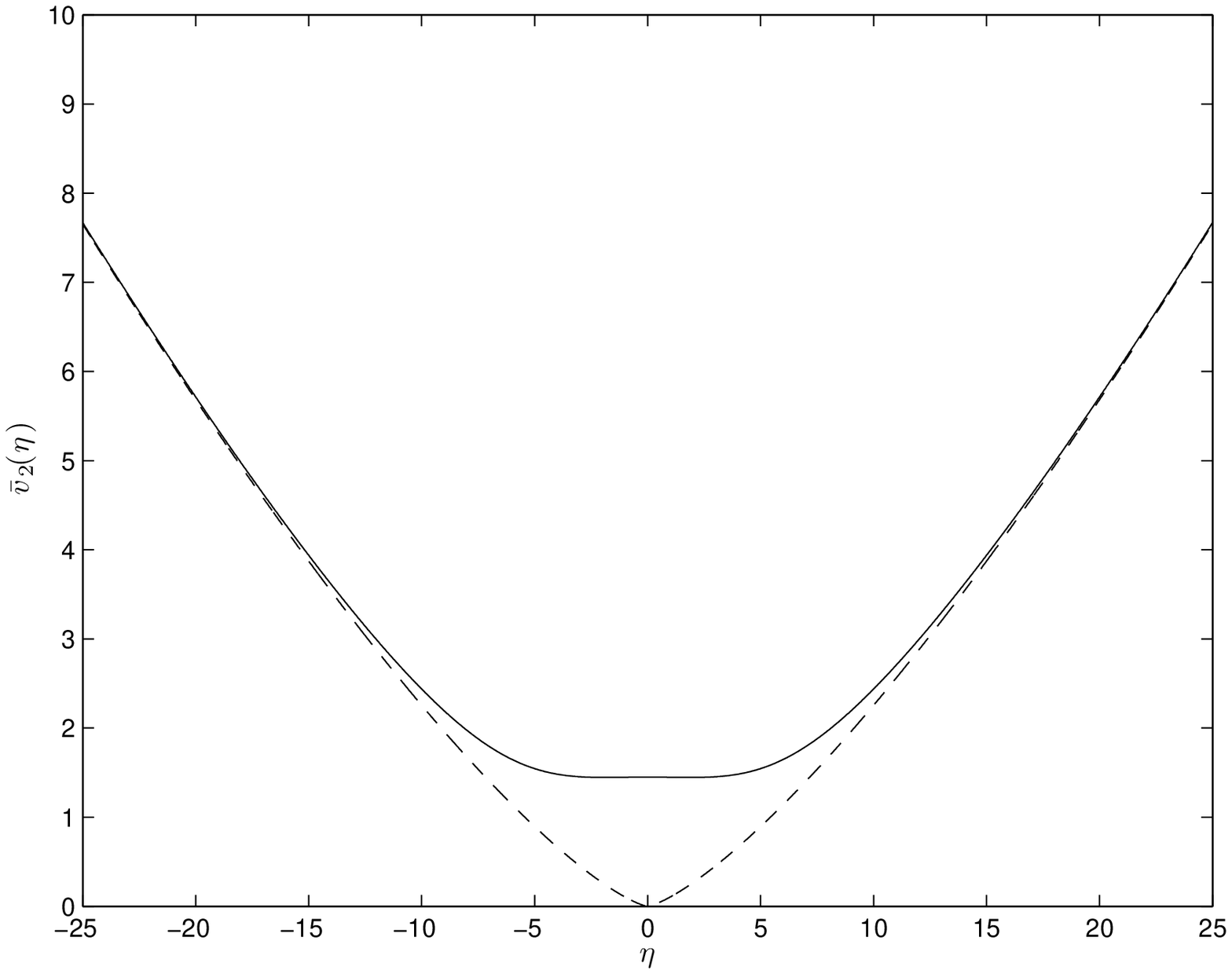}\label{fig_num_TDb}}
\parbox{4.5in}{\caption{Plots of two self-similar profiles $\bar{v}_1(\eta)$ and $\bar{v}_2(\eta)$ satisfying \eqref{main_v4}. The dotted curves represent the far field behaviour $\bar{v}_j(\eta)\sim c^{(j)}_0|\eta|^{4/3}$ as $|\eta|\to\infty$. The values $c^{(1)}_0 = 0.906$ and $c^{(2)}_0 = 0.1047$ were determined numerically. Note that $\bar{v}_2$ has a small dimple at the origin indicating three critical points.\label{fig_num_TD}}}
\end{figure}
After seeking convergence over a wide range of parameters $c_0$, exactly two solutions to \eqref{main_v4}, denoted $\bar{v}_1(\eta)$, $\bar{v}_2(\eta)$ were found \jl{as shown in} Fig.~\ref{fig_num_TD}. This solution multiplicity appears to be qualitatively similar in character to that observed \cite{BGW,GJFW,GALAK} in the self-similar blow up of fourth order PDEs with power law nonlinearity. To address the question of the existence of a stable self-similar quenching profile for \eqref{mems_main}, the linear stability of $\bar{v}_1(\eta)$ and $\bar{v}_2(\eta)$ is now analyzed by setting $ v = \bar{v}(\eta) + \phi(\eta)e^{\mu s}$ for $\phi\ll1$ in \eqref{main_v2} to arrive at the eigenvalue problem
\bsub\label{TD_eigenvalue}
\begin{align}
\label{TD_eigenvalue_a} \mu \phi &= -\phi_{\eta\eta\eta\eta} -\frac{\eta}{4}\phi_{\eta} + \left(\frac{1}{3} + \frac{2}{\bar{v}^3}\right)\phi, \qquad -\infty<\eta<\infty;\\[5pt]
\label{TD_eigenvalue_b} \mu \phi &= \frac{\phi}{3} -\frac{\eta}{4} \phi_{\eta}, \qquad \eta \to\pm \infty.
\end{align}
\esub
Apart from the following two modes associated with translation in touchdown time $t_c$ and location $x_c$
\begin{equation}\label{TD_eigenvalue_Spurious}
\mu_0 = 1, \qquad \phi_0 = \frac{\bar{v}}{3} -\frac{\eta}{4} \bar{v}_{\eta}; \qquad \mu_1 = \frac{1}{4}, \qquad \phi_1 =  \bar{v}_{\eta},
\end{equation}
the spectra of \eqref{TD_eigenvalue} must in general be determined numerically by reducing, via discretization, \eqref{TD_eigenvalue} to a linear system $\bar{\mathcal{L}}_{\mu} \phi=0$ and then seeking $\mu$ such that $\det\bar{\mathcal{L}}_{\mu} =0$. The eigenvalues appear to be purely real and the largest eight numerically obtained eigenvalues associated with each of the two profiles $\bar{v}_1(\eta)$ and $\bar{v}_2(\eta)$ are displayed in Table.~\ref{Tab_Evals_Strip}.  In the spectra associated with each profile, the two eigenvalues identified in \eqref{TD_eigenvalue_Spurious} are present. Ignoring these particular values, it is observed that the spectrum associated to $\bar{v}_1$ is strictly negative, while the spectrum associated with $\bar{v}_2$ contains two positive eigenvalues.
\begin{table}
\centering
\begin{tabular}{|c|c|c|c|c|c|c|c|c|}
\hline
{}&$\mu_0$&$\mu_1$&$\mu_2$&$\mu_3$&$\mu_4$&$\mu_5$&$\mu_6$&$\mu_7$\\
\hline
$\bar{v}_1$&$1.0003$&$0.2499$&$-0.1369$&$-0.4328$&$-0.6089$&$-0.8431$&$-1.1431$&$-1.4251$\\[5pt]
$\bar{v}_2$&$1.0000$&$0.7740$&$\phantom{-}0.5347$&$\phantom{-}0.2499$&$-0.0828$&$-0.4464$&$-0.8269$&$-1.2151$\\[5pt]
\hline
\end{tabular}
\vskip0.5cm
\parbox{4.5in}{\caption{The first eight numerically obtained eigenvalues of \eqref{TD_eigenvalue} for the two profiles  $\bar{v}_1(\eta)$ and  $\bar{v}_2(\eta)$. The value of $L=50$ and a uniform discretization with $N=1000$ grid points were used. \label{Tab_Evals_Strip} }}
\end{table}

\begin{figure}
\centering
\subfigure[$\eps=0.2$ with $x_c=0$.]{\includegraphics[width=0.4\textwidth]{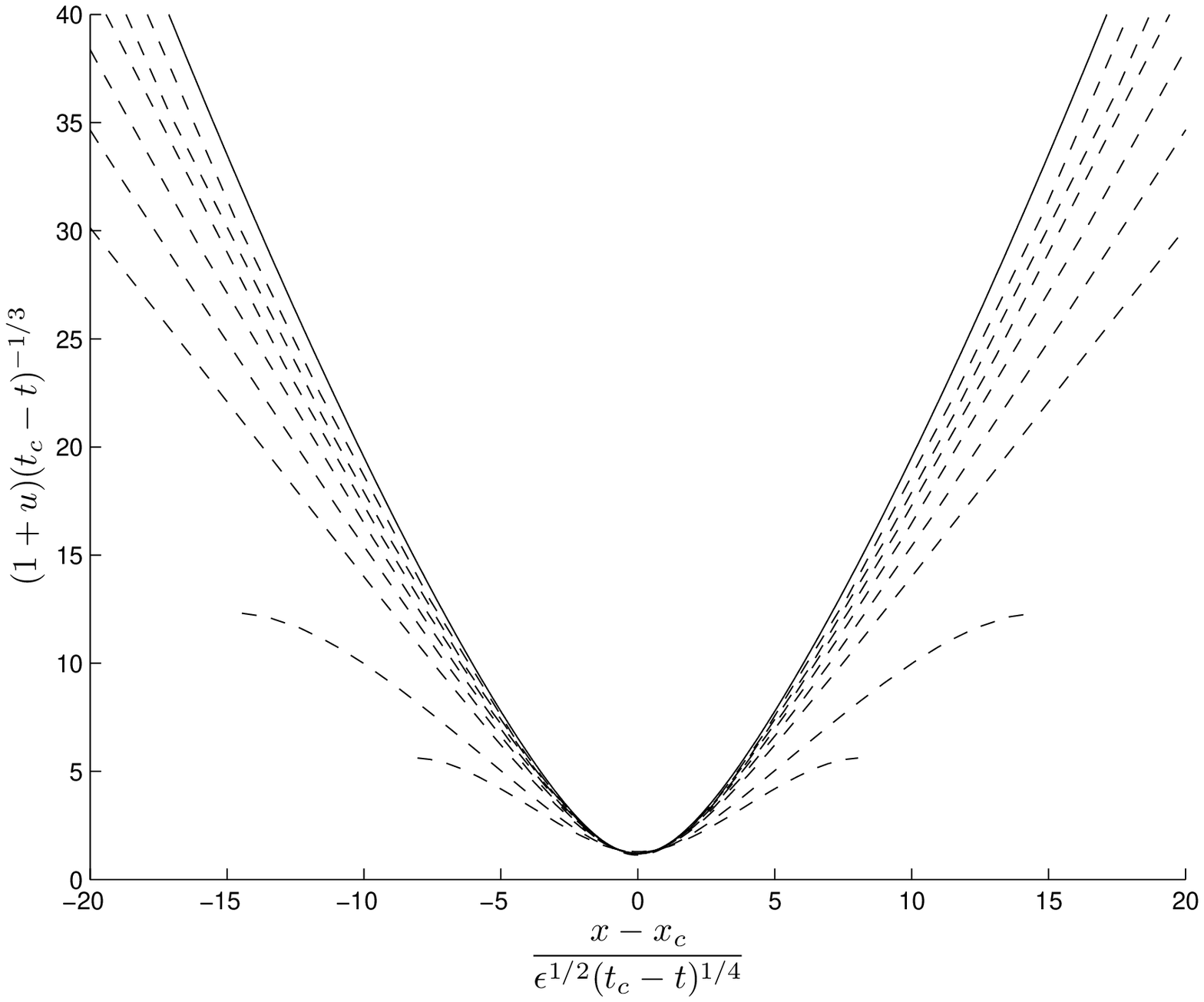}}\qquad
\subfigure[$\eps=0.02$]{\includegraphics[width=0.4\textwidth]{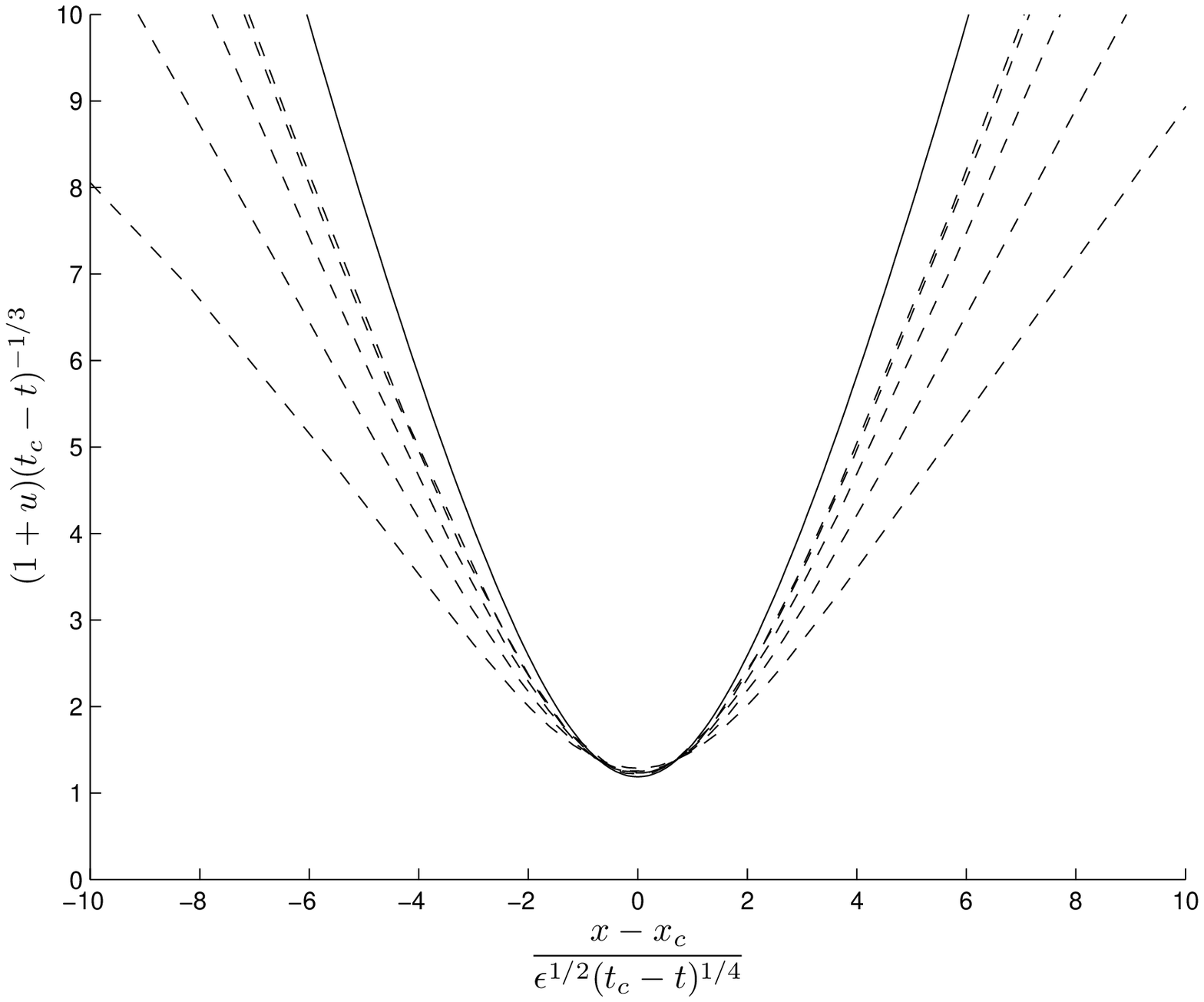}}
\parbox{4.5in}{\caption{Comparison of full numerical solutions (dashed lines) of \eqref{mems_main} to the stable self similar profile (solid line) $\bar{v}_1(\eta)$ for various $t$ approaching touchdown at $t_c$. This is for $N=1$ and clamped boundary conditions.\label{fig_num_CP} }}
\end{figure}

This suggests that the profile $\bar{v}_1(\eta)$ is a stable self-similar quenching profile for \eqref{mems_main} and indeed, in Fig.~\ref{fig_num_CP}, convergence of the full numerical solution to $\bar{v}_1(\eta)$ is observed as $t\to t_c^{-}$ for the case of touchdown at and away from the origin.

\subsection{Radially symmetric quenching solutions in 2D}

\begin{figure}
\centering
\subfigure[Profile $\bar{v}_1(\eta)$]{\includegraphics[width=0.4\textwidth]{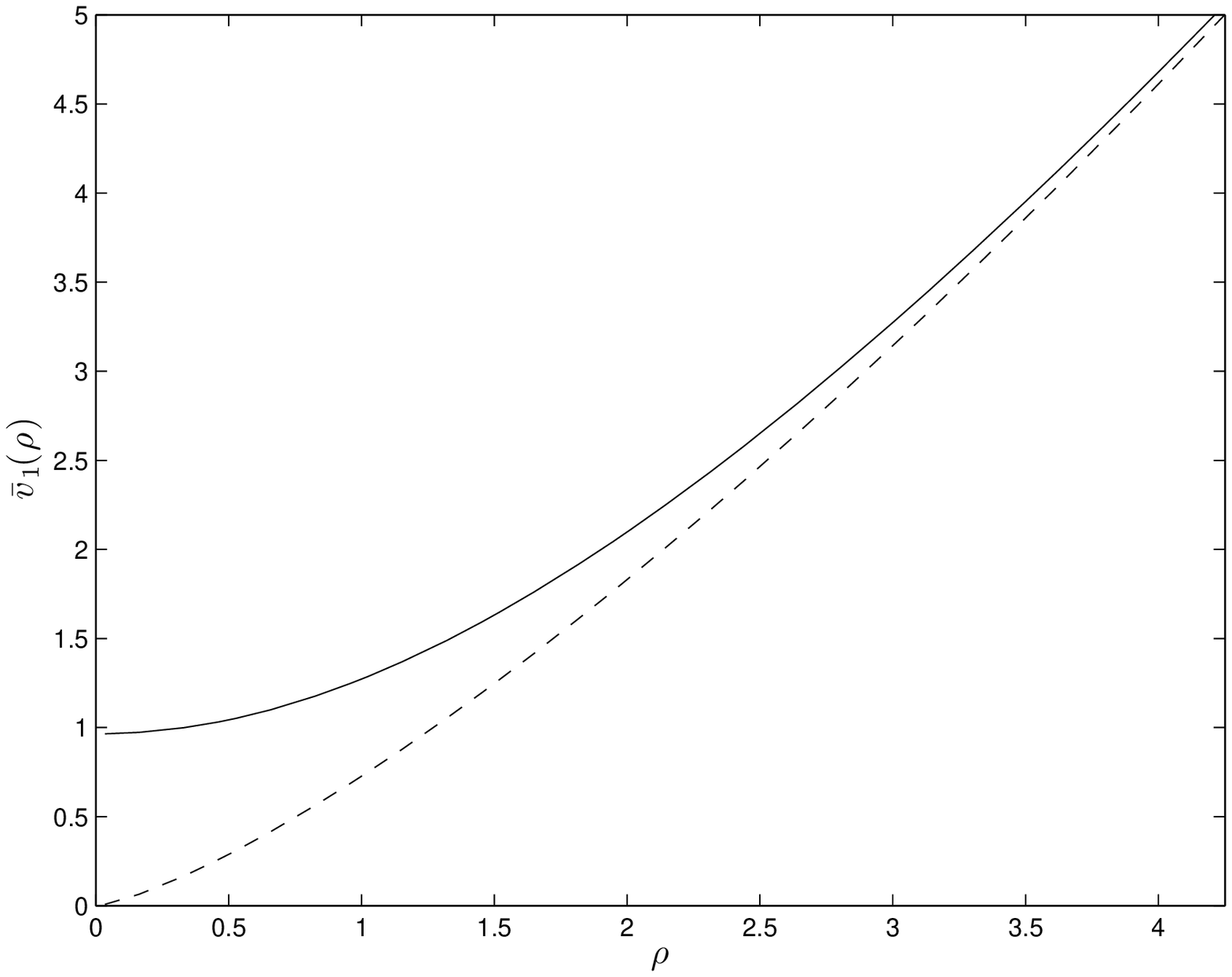}\label{radial_TDa}}\qquad
\subfigure[Profile $\bar{v}_2(\eta)$]{\includegraphics[width=0.4\textwidth]{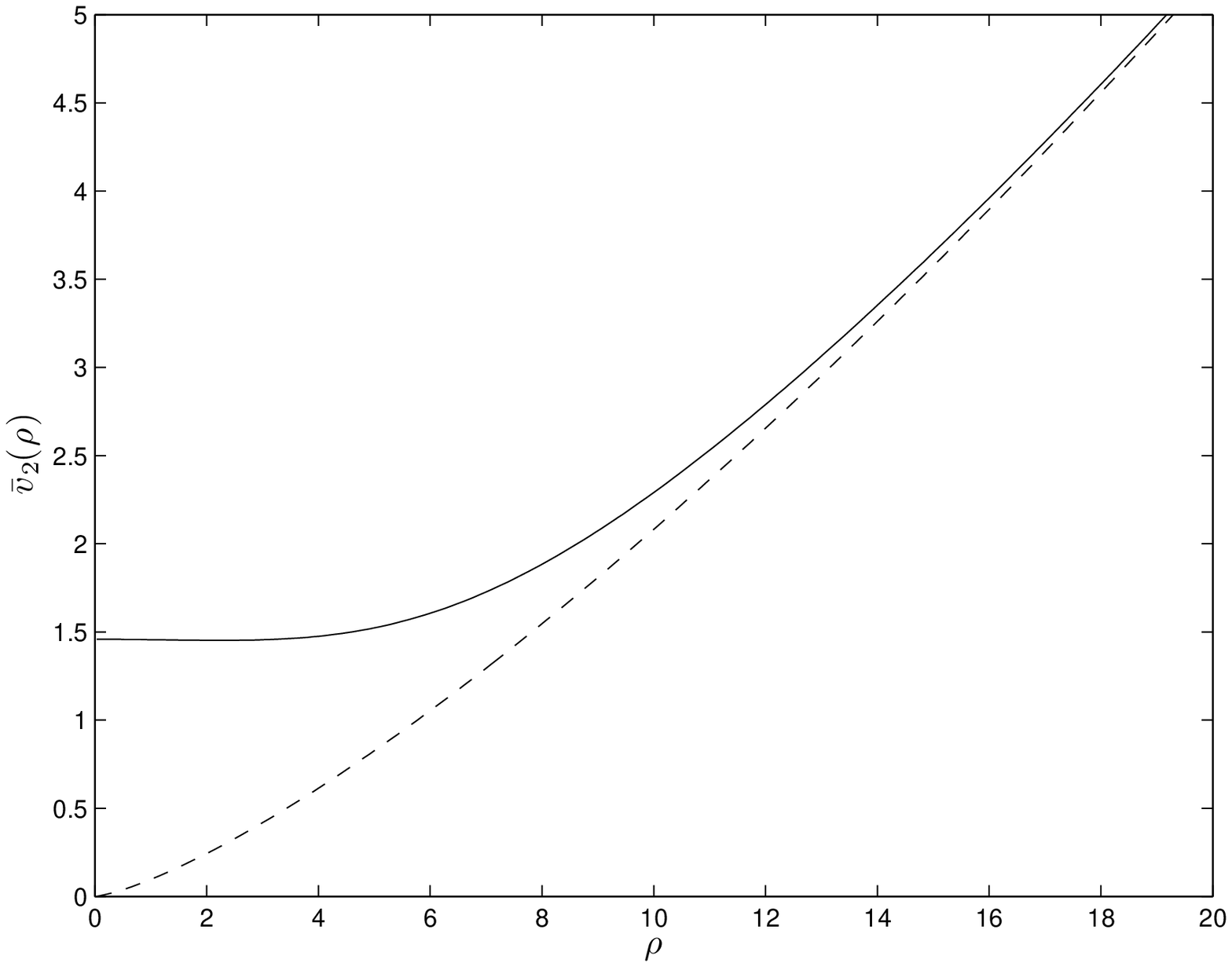}\label{radial_TDb}}
\parbox{4.5in}{\caption{Plots of two self-similar profiles $\bar{v}_1(\eta)$ and $\bar{v}_2(\rho)$ satisfying \eqref{Radial_eq5}. The dotted curves represent the asymptotic far field behaviour $\bar{v}_j(\rho)\sim c^{(j)}_0\rho^{4/3}$ as $\rho\to\infty$. The values $c^{(1)}_0 = 0.7265$ and $c^{(2)}_0 = 0.0966$ were determined numerically. Note that $\bar{v}_2$ has a small dimple at the origin indicating two critical points including that at $\rho=0$.\label{radial_TD}}}
\end{figure}

Self similar quenching profiles of the MEMS problem \eqref{mems_main} are now considered in two spatial dimensions. For radially symmetric solutions on the unit disc, the cases of touchdown at and away from the origin are treated separately. The variables
\begin{equation}\label{Radial_eq2}
u(x,t) = -1 + (t_c-t)^{1/3}v(\eta), \qquad \eta = \frac{r}{\eps^{1/2}(t_c-t)^{1/4}}, 
\end{equation}
which assume touchdown at the origin, transform \eqref{mems_main} to
\begin{equation}\label{Radial_eq3}
-\Delta_{\eta}^2 v -\frac{1}{4} \eta\cdot \nabla_{\eta} v + \frac{v}{3} - \frac{1}{v^2} = 0, \qquad \eta \in \mathbb{R}^2.
\end{equation}
which is a partial differential equation for the self similar quenching profile. The question of existence, multiplicity and stability  of solutions to \eqref{Radial_eq3} appears to be an open question in spatial dimensions $N\geq2$ (c.f. \cite{GALAK}). If a radially symmetric solution $v(\eta) = v(|\eta|)$ is presumed, then \eqref{Radial_eq3} reduces to 
\begin{equation}\label{Radial_eq4}
v'''' + \frac{2}{\rho}v''' - \frac{1}{\rho^2} v'' + \frac{1}{\rho^3}v' +\frac14 \rho\; v' - \frac{v}{3} +  \frac{1}{v^2} = 0
\end{equation}
where $\rho = |\eta|$.  A far field analysis similar to that which led to \eqref{main_v4} can be applied to \eqref{Radial_eq3} to establish boundary conditions which imply algebraic growth with exponentially growing terms suppressed at infinity. After algebra the full problem for the radially symmetric self-similar quenching profile in dimension $N=2$ is 
\bsub\label{Radial_eq5}
\begin{align}
\label{Radial_eq5_a}v'''' + \frac{2}{\rho}v''' -\frac{1}{\rho^2} v'' + \frac{1}{\rho^3}v' +\frac14 \rho\, v' - \frac{v}{3} +  \frac{1}{v^2} = 0 \qquad \rho>0\\[5pt]
\label{Radial_eq5_b}  v(\rho) \sim \left[ c_0 \rho^{4/3} + o(\rho^{4/3}) \right] + \bar{C}\rho^{-16/9}\exp\left[-3\rho^{4/3} 2^{-8/3} \right] + \cdots\qquad& \rho\to\infty
\end{align}
with symmetric conditions at the origin enforced by
\begin{equation}\label{Radial_eq5_c}
v'(0) = v'''(0) = 0.
\end{equation}
\esub
This nonlinear equation is solved numerically by first discretizing \eqref{Radial_eq5_a} on $[0,L]$ for $L$ large, applying far field behaviour \eqref{Radial_eq5_b} as a Robin condition ( c.f. \eqref{main_v3_b} ) at consecutive endpoints followed by Newton iterations initialized with \eqref{num_init}. The iterations are initialized over a wide range of the parameter $c_0$, with convergence observed for two isolated values $c_0^{(1)} = 0.7265$ and $c^{(2)}_0 = 0.0966$. The two associated profiles are displayed in Fig.~\ref{radial_TD}.

\begin{figure}
\centering
\subfigure[$\eps = 0.1$]{\includegraphics[width=0.4\textwidth]{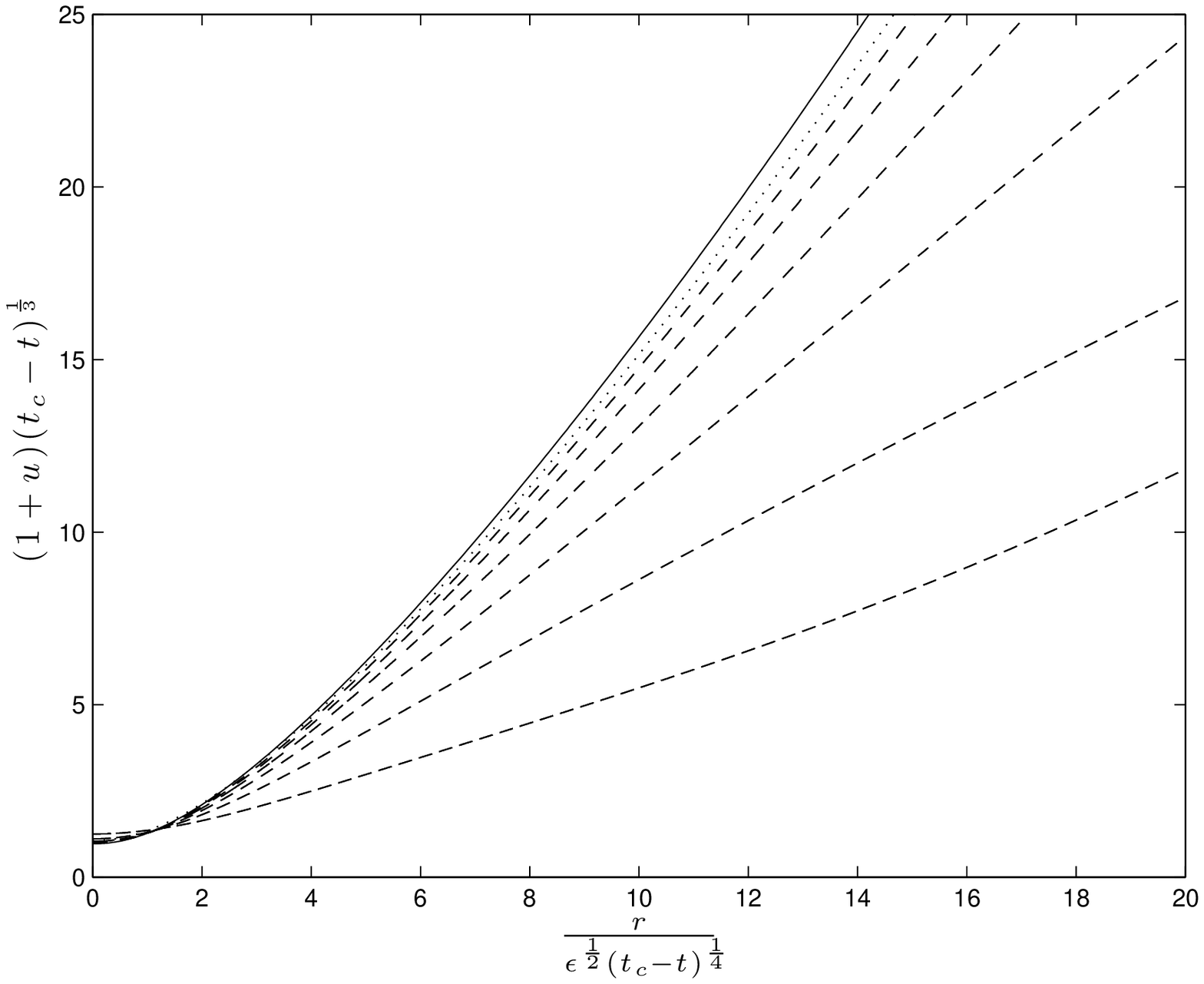}\label{radial_TD_comp_a}}\qquad
\subfigure[$\eps = 0.02$]{\includegraphics[width=0.4\textwidth]{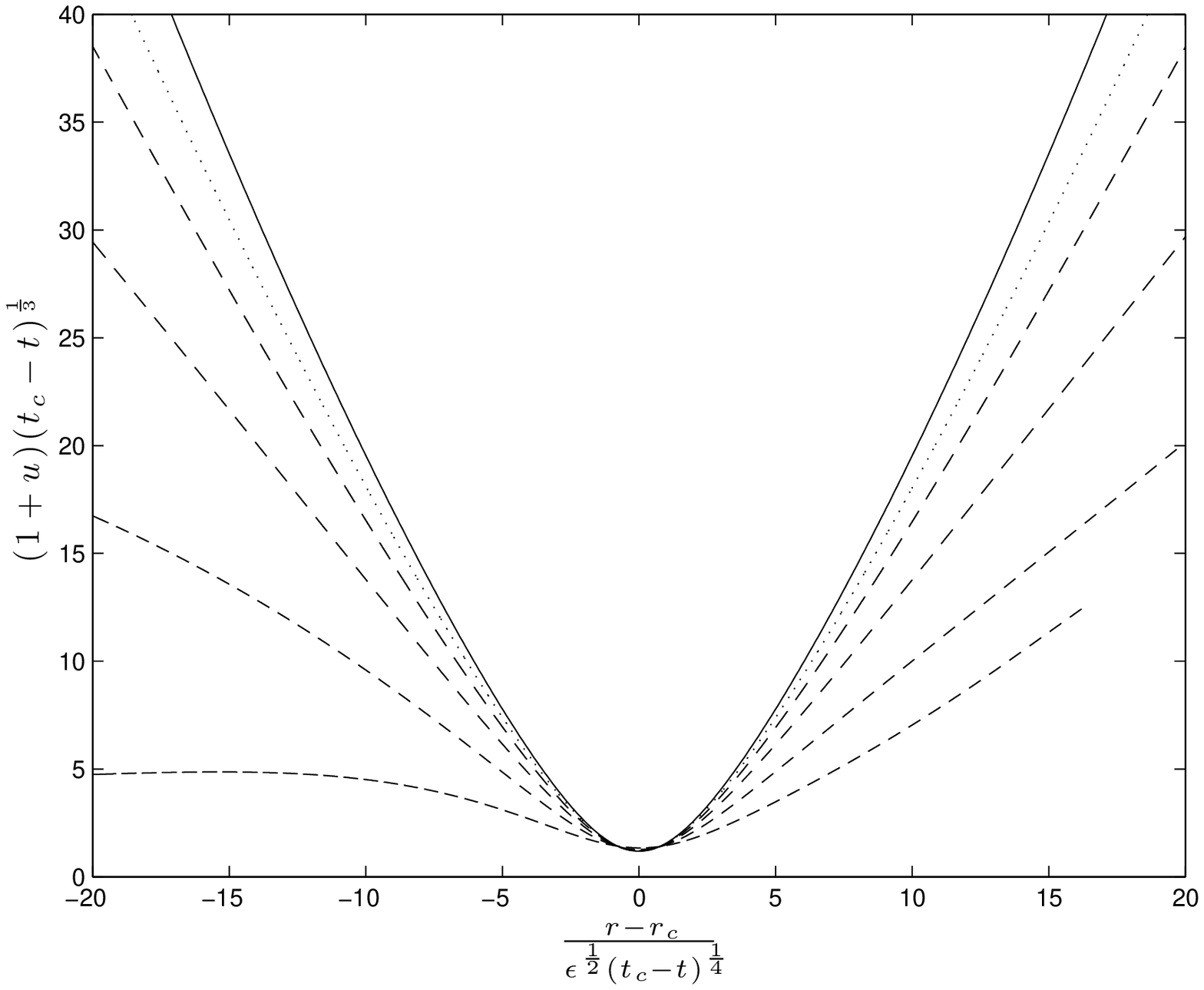}\label{radial_TD_comp_b}}
\parbox{4.5in}{\caption{Convergence of radially symmetric solutions (dashed \& dotted) of \eqref{mems_main} to self-similar profiles. Left: Touchdown is at the origin and convergence is observed to the monotone profile $\bar{v}_1(\eta)$ (solid) solving \eqref{Radial_eq5}. Right: For this case, touchdown is away from the origin and so convergence is to the monotone profile $\bar{v}_1(\eta)$ solving \eqref{main_v4}, as in the 1D strip case. In both figures, the dotted curve represents the solution for smallest $(t_c-t)$.}\label{radial_TD_comp}}
\end{figure}

As in the $N=1$ case, the second profile $\bar{v}_2(\eta)$ has a dimple at the origin and, as illustrated in Fig.~\ref{radial_TD_comp_a}, full numerical solutions of \eqref{mems_main} are observed to converge to the monotonic self-similar profile $\bar{v}_1(\eta)$. 

For touchdown away from the origin in the radially symmetric unit disc case, the self-similar quenching profile appears to be the same as that obtained for the $N=1$ case. Indeed, the appropriate similarity variables are 
\begin{equation}\label{Radial_eq6}
u(x,t) = -1 + (t_c-t)^{1/3}v(\eta), \qquad \eta = \frac{r-r_c}{\eps^{1/2}(t_c-t)^{1/4}}. 
\end{equation}
These variables rescale the biharmonic term as follows
\[
-\eps^2\Delta^2 u \to -(t_c-t)^{-2/3} \big[ v_{\eta\eta\eta\eta} + \mathcal{O}(\eps^{1/2}(t_c-t)^{1/4}) \big]
\] 
and so in the limit as $t\to t_c$, the $v_{\eta\eta\eta\eta}$ term is dominant. This results in a self-similar profile which satisfies 
\bsub\label{main_v4_radial}
\begin{align}
\label{main_v4_radial_a}  \bar{v}_{\eta\eta\eta\eta} +\frac{\eta}{4}\bar{v}_{\eta} - \frac{\bar{v}}{3} + \frac{1}{\bar{v}^2}&=0, \quad -\infty<\eta<\infty.\\[5pt]
\label{main_v4_radial_b}  \bar{v}\sim \sum_{n=0}^{\infty} c_n\eta^{4/3 -4n} &+ \bar{C}|\eta|^{-10/9}\exp\left[-3|\eta|^{4/3} 2^{-8/3} \right], \qquad \eta \to \pm\infty
\end{align}
\esub
as derived for the 1D case in \eqref{main_v4}. Consequently, quenching solutions away from the origin in the radially symmetric unit disc case are expected to converge to the self-similar quenching profile of the $N=1$ case, as  confirmed by the numerical simulations displayed in Fig.~\ref{radial_TD_comp_b}.

\section{Discussion}

Quenching solutions of a fourth order parabolic differential equation with a singular nonlinearity have been analyzed for a 1D strip and under radial symmetry on the unit disc with both clamped and Navier boundary conditions. In contrast to its second order equivalent, the fourth order PDE can quench at multiple points away from the origin. More precisely, in the case $N=1$, we have shown that the PDE can quench at two distinct points symmetric about the origin, while in the radially symmetric unit disc case, it can quench on an inner circle of finite radius.  In each case, the location of the quenching set was predicted by means of an asymptotic expansion whose accuracy is verified using very accurate adaptive numerical methods. 

In the limit as $t\to t_c$, where $t_c$ is the quenching time, the behaviour of \eqref{mems_main} was shown to be self-similar in nature. This is again in contrast with the second order equivalent of \eqref{mems_main}. The self-similar profile itself was obtained numerically and its limiting behaviour for $t\to t_c$ is given by
\begin{equation*}
u(x,t_c) = -1 + c_0\left(\frac{|x-x_c|}{\eps^{1/2}}\right)^{4/3}
\end{equation*}
where $c_0 = 0.7265$ for $N=2$ and touchdown at the origin and $c_0 = 0.9060$ in the other cases.

There are many interesting questions which stem from this study. In the case of the unit disc, it may be possible for the dynamics of \eqref{mems_main} to break the radial symmetry of the quenching set. All the simulations presented here were initialized with $u(x,0)=0$. Adding random noise to the initial condition breaks the left-right symmetry for $N=1$ and rotational symmetry for $N=2$. The symmetry breaking can be amplified by the dynamics of the PDE.

 In such a scenario, the ring would most likely be split up into a collection of points whose arrangement would need to be determined. The prediction of the quenching set of \eqref{mems_main} for larger classes of 2D geometries is another interesting open problem. For regular geometries, it may be that the number of axes of symmetry determines the quenching set but for irregular domains, it is not clear that the touchdown locations can be determined by simple geometric considerations. This question may be amenable to perturbation analysis, for example an almost circular domain whose boundary is $r= 1+ \delta f(\theta)$ for some $\delta\ll1$ and $f(\theta)$ a $2\pi$ periodic function. 

A robust method for solving \eqref{mems_main} of a large class of 2D geometries would be an essential complement to any analytical investigation of the above questions. In particular, a meshless method might be well suited to handle the highly non-uniform grids needed to resolve the dynamics of \eqref{mems_main} very close to touchdown \cite{MME}.

The treatment of these open issues is beyond the scope of this manuscript and will be left for future investigation.

\section{Acknowledgements}
A.E.L is very grateful to M.J. Ward for many useful discussions.

\begin{appendix}

\section{Spatial Discretization}\label{AppendixA}

For discretization in space, a collocation method based on piecewise 7th-order polynomial interpolation is employed. On the interval $x\in[X_i(t),X_{i+1}(t)]$ for $i = 0,1,\ldots,N$, the solution $u(x,t)$ is written as
\bsub\label{num_spatial}
\begin{equation}\label{num_spatial_a}
u(x,t) = \sum_{k=0}^{3} [u_{i}^{(k)}(t)L_{0,k}(s_i) + u_{i+1}^{(k)}(t) L_{1,k}(s_i) ] H_i^k
\end{equation}
where
\begin{equation}\label{num_spatial_b}
s_i =\frac{x-X_{i}(t)}{H_i(t)} \in [0,1], \qquad H_i(t) = X_{i+1}(t) - X_{i}(t),\qquad u_i^{(k)}(t) = \left[\frac{d^k}{dx^k} u(x,t)\right]_{x=X_i(t)} 
\end{equation}
and the $L_{0,j}(s)$, $L_{1,j}(s)$ for $j=0,1,2,3$ are the shape functions
\begin{equation}\label{num_spatial_c}
\begin{array}{ll}
L_{0,0}(s) = (20s^3 +10s^2+4s+1)(s-1)^4, & L_{0,1}(s) = s(10s^2+4s+1)(s-1)^4,\\[5pt]
L_{0,2}(s) = \ds\frac{s^2}{2}(4s+1)(s-1)^4, & L_{0,3}(s) = \ds\frac{s^3}{6}(s-1)^4,\\[5pt]
L_{1,0}(s) = -s^4(20s^3 -70s^2+84s-35), & L_{1,1}(s) = s^4(s-1)(10s^2-24s+15),\\[5pt]
L_{1,2}(s) = -\ds\frac{s^4}{2}(s-1)^2(4s-5) & L_{1,3}(s) = \ds\frac{s^4}{6}(s-1)^3.
\end{array}
\end{equation}
\jl{They satisfy}
\begin{equation}\label{num_spatial_d}
\left[\frac{d^p}{dx^p }L_{i,k}(s_i)\right]_{x=X_j(t)} = \quad \left\{ \begin{array}{cl} 1 & \mbox{if} \; i=j \; \mbox{and} \; k=p\\ 0 & \mbox{otherwise} \end{array},\right.
\end{equation}
\esub
\jl{so that the unknown coefficients $u_i^{(k)}(t)$ are the values of $u$ and its first three spatial derivatives at the nodal points $x = X_i(t)$. By construction, these are continuous at the nodal points.}

\jl{The dynamics of the $u_i^{(k)}(t)$ is obtained by substituting expansion \eqref{num_spatial_a} into the PDE, using the following expressions for the temporal and spatial derivatives of $u$.}
\bsub\label{num_spatial1}
\begin{align}
\label{num_spatial1_a}\frac{\partial^j}{\partial x^j} u(x,t) &= \sum_{k=0}^{3} \left[u_{i}^{(k)}(t)\frac{d^j}{ds^j} L_{0,k}(s_i) + u_{i+1}^{(k)}(t) \frac{d^j}{ds^j} L_{1,k}(s_i) \right] H_i^{k-j},\\
\nonumber \frac{\partial}{\partial t} u(x,t) &= \sum_{k=0}^{3}  \left[\frac{d}{dt}u_{i}^{(k)}(t)L_{0,k}(s_i) + \frac{d}{dt}u_{i+1}^{(k)}(t) L_{1,k}(s_i) \right] H^k_i\\
 \label{num_spatial1_b} &+ \frac{dH_i}{dt}\sum_{k=1}^{3}  \left[u_{i}^{(k)}(t)L_{0,k}(s_i) + u_{i+1}^{(k)}(t) L_{1,k}(s_i)\right]k H_i^{k-1} \\
\nonumber &- u_x(x,t)\left[\frac{dX_i}{dt} + s_i \frac{dH_i}{dt}\right].
\end{align}
\esub
Navier and clamped boundary conditions can be applied at both endpoints by choosing 
\begin{align}
\nonumber\mbox{(Clamped)}\quad u^{(0)}_0 = u^{(1)}_0 = u^{(0)}_{N+1} = u^{(1)}_{N+1}=0\\
\nonumber \mbox{(Navier)}\quad u^{(0)}_0 = u^{(2)}_0 = u^{(0)}_{N+1} = u^{(2)}_{N+1}=0
\end{align}
in the strip $\Omega = [-1,1]$ case and 
\begin{align}
\nonumber\mbox{(Clamped)}\quad u^{(1)}_0 = u^{(3)}_0 = u^{(0)}_{N+1} = u^{(1)}_{N+1}=0\\
\nonumber \mbox{(Navier)}\quad u^{(1)}_0 = u^{(3)}_0 = u^{(0)}_{N+1} = u^{(2)}_{N+1}+u^{(1)}_{N+1}=0
\end{align}
in the unit disc $\Omega = \{x^2+y^2\leq1\}$ case.

and the remaining \jl{equations} are \jl{obtained by writing the discretized PDE} at the Gauss points
\begin{equation}\nonumber
\rho_1 = \ds\frac{1}{2} - \frac{\sqrt{525+70\sqrt{30}}}{70}, \quad \rho_2 =  \ds\frac{1}{2} - \frac{\sqrt{525-70\sqrt{30}}}{70}, \quad \rho_3 = 1- \rho_2, \quad \rho_4= 1- \rho_1  
\end{equation}
on each interval $[X_i(t),X_{i+1}(t)]$, for $i=0,\ldots,N$. This provides $4(N+1)$ equations, which together with the four boundary conditions, \jl{are integrated in time to obtain} the $4(N+2)$ unknown nodal values $u_i^{(k)}(t)$ for $k=0,1,2,3$ and $i=0,\ldots,N+1$.
\end{appendix}

\end{document}